\documentclass[11pt]{article}

\usepackage{amsfonts}
\usepackage{fancyhdr}
\usepackage{comment}
\usepackage[a4paper, top=2.25cm, bottom=2.5cm, left=2.2cm, right=2.2cm]
{geometry}
\usepackage{times}
\usepackage{amsmath}
\usepackage{changepage}
\usepackage{amssymb}
\usepackage{graphicx}
\usepackage{mathrsfs}
\usepackage{relsize}
\usepackage{multicol}
\usepackage{enumerate}
\usepackage{tikz-cd}
\usepackage{color}
\usepackage[light,frenchstyle,narrowiints]{kpfonts}
\DeclareMathAlphabet\mathfrak{U}{euf}{m}{n}
\SetMathAlphabet\mathfrak{bold}{U}{euf}{b}{n}
\usepackage[framemethod=tikz]{mdframed}
\usepackage{amsthm}
\usepackage{titlesec}
\usepackage[rightcaption]{sidecap}
\sidecaptionvpos{figure}{c}
\usepackage{scrextend}
\usepackage{comment}

\theoremstyle{definition}
\newmdtheoremenv[
  hidealllines=true,
  innerleftmargin=10pt,
  innerrightmargin=10pt,
  innertopmargin=0pt,
]{defi}{Definition}
\theoremstyle{plain}
\newmdtheoremenv[
  hidealllines=true,
  innerleftmargin=10pt,
  innerrightmargin=10pt,
  innertopmargin=0pt,
]{thm}[defi]{Theorem}
\theoremstyle{plain}
\newmdtheoremenv[
  hidealllines=true,
  innerleftmargin=10pt,
  innerrightmargin=10pt,
  innertopmargin=0pt,
]{lemma}[defi]{Lemma}
\theoremstyle{plain}
\newmdtheoremenv[
  hidealllines=true,
  innerleftmargin=10pt,
  innerrightmargin=10pt,
  innertopmargin=0pt,
]{cor}[defi]{Corollary}
\theoremstyle{definition}
\newmdtheoremenv[
  hidealllines=true,
  innerleftmargin=10pt,
  innerrightmargin=10pt,
  innertopmargin=0pt,
]{ex}[defi]{Example}

\makeatletter
\renewcommand*\env@matrix[1][*\c@MaxMatrixCols c]{%
  \hskip -\arraycolsep
  \let\@ifnextchar\new@ifnextchar
  \array{#1}}
\makeatother

\newcommand{\R}{\mathbb{R}}

\newcommand{\B}{\mathbb{B}}

\newcommand{\bs}{\bigskip}
\newcommand{\bbs}{\bigskip\bigskip}

\newcommand{\ml}[1]{\mathlarger{#1}}

\newcommand{\mbb}[1]{\mathbb{#1}}

\newcommand{\mcal}[1]{\mathcal{#1}}
\newcommand{\txt}[1]{\textrm{#1}}
\newcommand{\cb}[1]{ \big \{ #1 \big \} }
\newcommand{\cp}[1]{ \big ( #1 \big ) }
\renewcommand{\inf}[1]{\textrm{inf}\big \{ #1 \big \}}
\renewcommand{\sup}[1]{\textrm{sup}\big \{ \es #1 \big \}}
\newcommand{\es}{\enskip}

\newcommand*{\LargerCdot}{\raisebox{0ex}{\scalebox{0.45}{$\bullet$}}}
\newcommand{\ba}[2]{\left [\begin{array}{#1} #2 \end{array}\right ]}
\renewcommand{\dot}[1]{\overset{\LargerCdot}{#1}}
\titleformat{\section}[block]
{\normalfont\sffamily}
{\thesection}{.5em}{\titlerule\\[.8ex]\bfseries}
\titleformat{\subsection}
  {\normalfont\scshape}{\thesubsection}{1em}{\bfseries}
\newcommand\restr[2]{{
  \left.\kern-\nulldelimiterspace
  #1
  \vphantom{\big|}
  \right|_{#2}
  }}

\begin{document}


\title{\vspace*{-1.55cm}\textothersc{A Variational Approach\\ to Local Asymptotic and Exponential Stabilization of Nonlinear Systems}}
\author{\small{ BRYCE A. CHRISTOPHERSON$\ml{\vphantom{\big|}^*}$\let\thefootnote\relax\footnote{Corresponding author \newline \hspace*{0.25cm} 2010 Mathematics Subject Classification: 93D15, 93CF10, 49J53 \newline \hspace*{0.25cm} Key Words: nonlinear control systems, stabilization by feedback, asymptotic and exponential stabilizability, variational analysis, linear openness}, \es  FARHAD JAFARI, \es BORIS S. MORDUKHOVICH\footnote{Research of this author was partly supported by the USA National Science Foundation under grants DMS-1512846 and DMS-1808978, by the USA Air Force Office of Scientific Research under grant \#15RT04, and by the Australian Research Council under Discovery Project DP-190100555.}}}
\date{}
\maketitle


\vspace*{-1.5cm}\section*{\textothersc{Abstract}}

Local asymptotic stabilizability is a topic of great theoretical interest and practical importance. Broadly, if a system $\dot{x} = f(x,u)$ is locally asymptotically stabilizable, we are guaranteed a feedback controller $u(x)$ that forces convergence to an equilibrium for trajectories initialized sufficiently close to it.  A necessary condition was given by Brockett:  such controllers exist only when $f$ is open at the equilibrium. Recently, Gupta, Jafari, Kipka and Mordukhovich \cite{GJKM} considered a modification to this condition, replacing Brockett's topological openness by the {\em linear openness} property of modern {\em variational analysis}. In this paper, we show that under the linear openness assumption there is a local diffeomorphism of neighborhoods of the equilibirum on which the system is exponentially stabilizable by means of continuous stationary feedback laws.  Introducing a transversality property and relating it to the above diffeomorphism, we prove that linear openness and transversality on a {\em punctured neighborhood of an equilibrium} is sufficient for local exponential stabilizability of systems with a rank deficient linearization.The main result goes beyond the usual Kalman and Hautus criteria for the existence of exponential stabilizing feedback laws, since it allows us to handle systems for which  exponential stabilization is achieved through higher-order terms. However, it is implemented so far under a rather restrictive row-rank conditions. This suggests a twofold approach to the use of these properties:  a {\em point-wise} version is enough to ensure stability via the linearization, while a {\em local} version is enough to overcome deficiencies in the linearization.


\vspace*{-0.5cm}\section{Introduction}


Consider the autonomous control system governed by nonlinear ordinary differential equations

\vspace*{-0.5cm}\begin{align}
    \dot{x}=f(x,u), & \es t \geq 0, \label{sys}
\end{align}

\vspace*{-0.2cm}\noindent where $f : \mcal{U} \times \mcal{V} \rightarrow \R^n$ on the right-hand side of \eqref{sys} is sufficiently smooth on some open set \newline $\mcal{U} \times \mcal{V} \subseteq \R^n \times \R^m$ containing the origin as a given equilibrium pair. The continuous feedback stabilizability property of our interest here is formulated as follows; see, e.g., \cite[Definition 10.11]{coron1}.

\begin{defi} \label{las}

    \noindent We say that the control system \eqref{sys} is \textit{locally asymptotically stabilizable} by means of {\it continuous stationary feedback laws} if there exists a continuous state-feedback control $\bar u(x)$ such that $\bar u(0)=0$ and the closed-loop system $\dot{x}=f(x,\bar u(x))$ has a locally asymptotically stable equilibrium at the origin.

\end{defi}


A variety of conditions describing whether system \eqref{sys} can be locally asymptotically stabilized by means of continuous stationary feedback laws have been derived; see, e.g., \cite{artstein,brockett,byrnes,coron2,coron1,hermes,sontag1,sontag2,sontag3}.  However, in the vast majority of cases the obtained conclusions are exclusively either necessary or sufficient, while those which do achieve both necessity and sufficiency occur only in more restrictive contexts or under more stringent conditions than what is posed in \eqref{sys}.  The most well-known result is a necessary condition given by Brockett \cite{brockett}, which ensures by using topological degree theory that such controllers exist only when $f$ satisfies a certain openness property.

\begin{thm} \label{brocketts}

    \noindent If system \eqref{sys} is locally asymptotically stabilizable by means of continuous stationary feedback laws, then it is necessary that $0 \in{\rm int}\,f\cp{\mcal{O}}$ for any neighborhood $\mcal{O} \subseteq \mcal{U} \times \mcal{V}$ of the origin, where `{\rm int}' signifies the set interior.
\end{thm}


In other words, any system that is asymptotically stabilizable by means of continuous stationary feedback laws is locally surjective onto a neighborhood of the origin.  The reader may note that, as in the inverse function theorem or implicit function theorem, it is clearly \textit{sufficient} for the Jacobian of $f$ to be of maximal rank to achieve this topological openness property, but it is certainly not {\it necessary}; take, e.g., $f(x)=x^3$. Having in mind this distinction, a partial converse in Brockett's theorem has been recently obtained by Gupta, Jafari, Kipka and Mordukhovich \cite{GJKM}. Their results emerge via a novel use of quantitative counterpart of the topological openness property that has been well recognized in variational analysis under the name of {\it linear openness} or, equivalently, {\it metric regularity} \cite{mordukhovich2,rockafellar}. However, this converse applies only to locally {\it exponentially stabilizable systems} understood in the following sense.
\begin{defi}\label{les}
    \noindent The control system \eqref{sys} is said to be \textit{locally exponentially stabilizable by means of continuous stationary feedback laws} if there exists a continuous stationary state-feedback control $\bar u(x)$ and constants $\alpha > 0$, $M> 0$, and $\delta > 0$ such that, for any starting point $x_0 \in \R^n$ with $\| x_0 \| < \delta$ there exists a unique solution $\phi(t, x_0)$ to the closed-loop control system $\dot{x}=f(x,\bar u(x))$ satisfying $\phi(0, x_0) = x_0$ and the exponential decay condition
    
    $$\|\phi(t, x_0)\| \leq M\|x_0\|e^{-\alpha t} \es \textrm{whenever} \es t \geq 0.$$
\end{defi}

A major goal of this paper is to develop a {\it variational approach} to investigate local asymptotic and exponential stabilizability of nonlinear control systems. In this way we show, in particular,  that systems which are {\it linearly open}--and hence satisfy, due to Theorem~\ref{brocketts}, a necessary condition for local asymptotic stabilizability by means of continuous stationary feedback laws--are locally exponentially stabilizable under the action of a {\it composition operator} generated by a {\it local diffeomorphism}. These developments go beyond the classical Kalman and Haustus criteria for exponential feedback stabilizability by allowing us  to achieve exponential and asymptotic stabilization by using higher-order terms. However, the realization of the variational approach in this paper requires imposing a restrictive row-rank assumption, which we are going to overcome in future research.

One may note that the application of a feedback law $u(x)$ to the system \eqref{sys} may be viewed as the action of a composition operator $f(x,u) \mapsto f(h(x,u))$ with $h(x,u)=(x,u(x))$.  As such, this is a natural (albeit much weaker) generalization of the concept of feedback stabilizability.  While it is clear that the existence of a composition operator stabilizing the system does not, in general, guarantee the stabilizability of the system by means of feedback laws, we will show that the imposition of certain additional conditions render the two equivalent. Namely, we introduce several forms of a particular {\it transversality condition} under which the linear openness proves to be {\it sufficient} for both asymptotic and exponential stabilizability.  This \textit{transversality condition} imposes that the geometry of the tangent bundle must be of an amicable form, which in turn reveals that the local geometry of the problem is sufficiently agreeable to allow the much weaker stabilizability-by-composition property to agree with the more strict continuous feedback stabilizability property of Definition~\ref{las}.  This suggests that systems satisfying both the \textit{variational} linear openness extension of Brockett's theorem and the aforementioned diffeo-geometric transversality condition have a {\it preferred geometry}. Translation of the system to such a geometry can thus ensure asymptotic and exponential stabilizability and control in this preferred setting.\vspace*{0.05in}

The rest of the paper is organized as follows. Section~\ref{sec2} collects some fundamental notions and results of variational analysis that lay at the heart of our variational approach to stabilizability of nonlinear systems. Section~\ref{sec3} is mostly devoted to the study of sufficient conditions for local exponential stabilizability by using composition feedback operators. Section~\ref{sec4} is the culmination of the paper, which develops a partial characterizations of both exponential and asymptotic stabilizability of linearly open control systems by means of continuous stationary feedback laws under the various novel transversality conditions introduced therein.  Finally, in the conclusion of Section~\ref{sec4} these techniques are employed to characterize the stabilizability of certain systems with a \textit{rank-deficient Jacobian at the origin}, provided they satisfy certain linear openness and transversality conditions locally in a nearby region.  This final result, when coupled with the preceding, yields a twofold approach to characterizing local stabilizability in terms of linear openness and transversality conditions:  a \textit{pointwise} version \textit{at an equilibrium point} is enough to ensure local exponential stability via the linearization, while a local version \textit{near an equilibrium point} is enough to overcome whatever deficiencies the linearization may present.  Besides theoretical assertions, we present a number of numerical examples that illustrate the major features of the obtained results. The concluding Section~\ref{sec5} summarizes the main developments of the paper and discusses some open problems of future research.\vspace*{0.05in}

Throughout the paper we use standard notation of control theory and variational analysis; see. e.g., \cite{coron1,mordukhovich3,rockafellar}. Recall that $\B_X$ and $\mbb{S}_X$ denote, respectively, the closed unit ball and the unit sphere of the space $X$, $\B_r(x)$ stands for the closed ball centered at $x$ with radius $r>0$, $\mbb{P}_V$ denotes the orthogonal projection onto a closed linear subspace $V$ of $X$ (provided $X$ is an inner product space), and $A^*$ indicates the adjoint operator for $A$, or the transposition of the matrix $A$ in the case $A$ is real-valued. We always use the convention that $\sup{\emptyset} = -\infty$.


\section{Linear Openness and Related Properties of Nonlinear Mappings}\label{sec2}


We begin this section with recalling the topological openness property used in Brockett's theorem and then proceed with quantitative tools and results of variational analysis employed in what follows. The given definitions and facts are restricted to finite-dimensional spaces, which is the framework of this paper, while most  hold in infinite dimensions as well; see, e.g., \cite{mordukhovich2}.

\begin{defi}\label{open} 
    \noindent A mapping $f:\R^\ell\rightarrow \R^n$ is said to be \textit{open at a point} $\bar{z}\in\R^\ell$ if $f(\bar{z}) \in \txt{int}f(\mcal{O})$ for any neighborhood $\mcal{O}$ of $\bar{z}$.
\end{defi}

As it is well known, the origin of this property goes back to the classical Banach-Schauder open mapping theorem.  It should be noted that \textit{open mappings} is a more restrictive class of mappings than those that are open about a point. The openness at a point property was used by Brockett \cite{brockett} via topological degree theory as stated in Theorem~\ref{brocketts}.  However, while a mapping may have the openness property at a given point, this property alone does not admit any measuring device to quantify the rate at which the mapping  is open about that point. The device that is most commonly utilized to achieve this refinement is \textit{linear openness} (known also as {\it covering} and as {\it openness at a linear rate}), which has been recognized as a fundamental property in variational analysis for nonlinear (as well as for non-smooth and set-valued) mappings; see, e.g., \cite{mordukhovich3,rockafellar}. This property is formulated as follows.

\begin{defi} \label{linear open}
    \noindent A mapping $f: \R^\ell \rightarrow \R^n$ is said to be \textit{linearly open} around $\bar{z}$ with modulus $\kappa > 0$ if there exists a neighborhood $\mcal{O}$ of $\bar{z}$ such that
    $$
    \B_{\kappa r}\big(f(z)\big)\subseteq f\big(\B_r(z)\big) \quad \txt{for any} \es z \in \mcal{O} \es \txt{and} \es r>0 \es \txt{with} \es \B_r(z) \subseteq \mcal{O}.
    $$
    \noindent Furthermore, we define the supremum over all the moduli $\kappa$ for which $f$ is linearly open around $\bar{z}$ to be the \textit{linear openness/covering bound} of $f$ around $\bar{z}$ and denote it by $cov f(\bar{z})$. If $f$ is not linearly open around $\bar{z}$ for any modulus $\kappa $, we put by convention $cov f(\bar{z}):=0$.
\end{defi}


Linear openness first appeared in \cite{milyutin} under the name of ``covering in a neighborhood" while it has been introduced and popularized by Milyutin in his talks and personal communications long before the publication of \cite{milyutin}. The interest to designating and calculating the (quantitative) exact bound $cov f(\bar{z})$ arose later on; see \cite{mordukhovich2,rockafellar,ioffe} and the references therein.

It is clear that the linear openness property from Definition \ref{linear open} yields its openness counterpart from Definition~\ref{open}, while the opposite implication fails. Let us emphasize that, in contrast to the openness property, the linear openness from Definition \ref{linear open}  ensures the {\it uniformity} of covering near $\bar z$ and designates a {\it linear rate} of change in the uniformity quantified by the modulus $\kappa$.

Linear openness is closely related to another measuring apparatus in variational analysis known as \textit{metric regularity}.  It has been well recognized that this fundamental property is in fact equivalent to linear openness with the reciprocal relationship between the exact bounds of the corresponding moduli; see Lemma~\ref{lomr identity} below for the precise formulation.

\begin{defi} \label{metric reg}
    \noindent A mapping $f: \R^\ell \rightarrow \R^n$ is said to be \textit{metrically regular} around $\bar{z}$ with modulus $\mu > 0 $ if there exists neighborhoods $\mcal{O}_1$ of $\bar{z}$ and $\mcal{O}_2$ of $\bar{y}:=f(\bar{z})$ such that we have the distance estimate
    $$
    d\big(z;f^{-1}(y)\big) \leq \mu \|y-f(z)\| \quad \txt{for any} \es z \in \mcal{O}_1 \es \txt{and} \es y \in \mcal{O}_2.
    $$
    where $d(z;\Omega) := \inf{ \|z-u\|\;\big|\; u \in \Omega}$ for a given set $\Omega$.
    \noindent Furthermore, we define the infimum over all the moduli $\mu$ for which $f$ is metrically regular around $\bar{z}$ to be the \textit{exact regularity bound} of $f$ around $\bar{z}$ and denote it by $reg f(\bar{z})$.  If $f$ is not metrically regular around $\bar{z}$ for any modulus $\mu$, we put by convention $reg f(\bar{z}): = \infty$.

\end{defi}


The interconnection between the two concepts above can be summarized by the following lemma; see, e.g., \cite[Theorem~3.2]{mordukhovich3} for a detailed proof of this and other relationships.

\begin{lemma} \label{lomr identity}  
    Assume that $cov f(\bar{z}) = \kappa > 0$ for an arbitrary mapping $f: \R^\ell \rightarrow \R^n$. Then we have that $reg f(\bar{z}) = \kappa^{-1} > 0$ and that $cov f(\bar{z}) \cdot reg f(\bar{z})=1$.
\end{lemma}


Variational analysis has achieved a complete characterization of linear openness and metric regularity (as well as the equivalent Lipschitz-like property of the inverses) for general set-valued mappings, with precise formulas for calculating the exact bounds of the corresponding moduli. This result is known in variational analysis and numerous applications as the co-derivative or Mordukhovich criterion for which the available proofs are based on variational/extremal principles of variational analysis; see \cite{mordukhovich1,mordukhovich3,rockafellar} and the references therein. In deriving the main results of the paper addressing smooth control systems, we use the following consequence of this criterion for $C^1$ mappings between finite-dimensional spaces and refer the reader to \cite[Theorem~1.57]{mordukhovich2} for an alternative device (in the general Banach space setting), where the proof of sufficiency is based on the Lyusternik-Graves iterative process.
\begin{thm} \label{lo fr jac}
    \noindent Let $f$ be of class $C^1$ in a neighborhood $\mcal{U} \times \mcal{V} \subseteq \R^\ell \times \R^n$ of the origin. Then $f$ enjoys the equivalent linear openness and metric regularity properties around the origin if and only if the Jacobian $\restr{J_f}{(0,0)}$ is of full rank.
    \noindent Furthermore, the exact bounds of the linear openness and metric regularity of $f$ at the origin are precisely given by, respectively,
    $$
    cov f(0,0)= \mathrm{min}\cb{\big\|\restr{J_f}{(0,0)}^*v\big\|\;\big|\;\|v\|=1}\quad{\rm and}\quad reg f(0,0) =\big\|\big(\restr{J_f}{(0,0)}^*\big)^{-1}\big\|.$$
\end{thm}


Let us also recall yet another definition from variational analysis that allows us to reformulate the properties given in Definitions~\ref{linear open} and  \ref{metric reg}; see \cite{ioffe} for more references and discussions. Although the definition below makes sense for operators between arbitrary Banach spaces, we formulate it just in the finite dimensional case that is needed for this paper.

\begin{defi}\label{b constant}
    \noindent Let $T:\R^\ell \rightarrow \R^n$ be a bounded linear operator. Then the \textit{Banach constant} of $T$ is defined as the quantity
    $$\Gamma(T):=\mathrm{sup}\big\{r\ge 0\;\big|\;r\B_{\R^n}\subseteq T\big(\B_{\R^l}\big)\big \}=\mathrm{inf}\big\{\|y\|\;\big|\;y\notin T\big(\B_{\R^\ell}\big)\big\}.
    $$
\end{defi}
\vspace*{0.2in}

\vspace*{-0.15cm}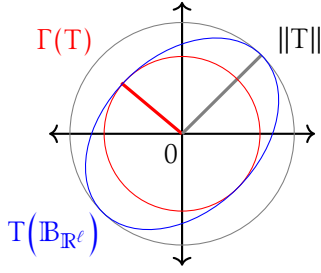
\begin{SCfigure}[2.5][ht]
\caption{The Banach constant $\Gamma(T)$ and norm $\|T\|$ of a linear operator $T$ shown as the radii of the balls contained in and containing the $T$-image of the unit ball $T(\B_{\R^\ell})$, respectively.} \begin{tikzpicture}[scale=0.5]
\draw (1.7,2) node[below]{$0$};
\draw[text=red] (-1.1,5.1) node[below]{$\Gamma(T)$};
\draw (5.1,5.1) node[below]{$\|T\|$};
\draw[text=blue] (-1.4,0.2) node[below]{$T\big (\mbb{B}_{\R^\ell} \big )$};
\draw[<->, thick] (2,-1.5) -- (2,5.5);
\draw[<->, thick] (-1.5,2) -- (5.5,2);
\draw[very thick,red] (2,2) -- (0.4,3.35); 
\draw[very thick,gray] (2,2) -- (4.09,4.09);
\draw[gray] (2,2) circle (2.95cm);
\draw[red] (2,2) circle (2.05cm);
\draw[blue,rotate around={-45:(2,2)}] (2,2) ellipse (2.05cm and 2.95cm);
\end{tikzpicture}\label{b constant fig} \end{SCfigure}

The Banach constant of a linear operator can be computed straightforwardly; see \cite{ioffe}.
\begin{lemma} \label{b constant formula} 
    Having $\|T^{-1}\|:=\mathrm{sup}\Big\{\mathrm{inf}\big\{\|x\|\;\big|\;T x = y\big\}\;\Big|\;\|y\|=1\Big\}$ for a bounded linear operator $T$, we get
    $$\Gamma(T)=\mathrm{inf}\Big\{\|T^*y\|\;\Big|\;\|y\|=1\Big\}=\frac{1}{\|T^{-1}\|}.$$
\end{lemma}


Combining Theorem~\ref{lo fr jac}, Lemma~\ref{b constant formula} with  Definitions~\ref{linear open}, \ref{metric reg}, and \ref{b constant} allows us to arrive at the following expressions for the exact bounds of linear openness and metric regularity of smooth mappings via the Banach constant.

\begin{cor}\label{cov b constant}
    \noindent Let $f: \R^\ell \rightarrow \R^n$ be of class $C^1$ in a neighborhood of $\bar{z} \in \R^l$. Then we have the exact bound formulas
    $$
    cov f(\bar{z}) = \Gamma\left ( \restr{J_f}{\bar{z}} \right ) \es {\rm and} \es reg f(\bar{z}) = \frac{1}{\Gamma \left (\restr{J_f}{\bar{z}} \right )} = \left \| \restr{J_f}{\bar{z}}^{-1} \right \|.$$

\end{cor}


This corollary shows that the exact linear openness bound of a smooth mapping $f$ is the Banach constant of its linearization, and that the exact regularity bound of $f$ is the operator norm of the inverse of the linearization. Overall, the variational results presented in this section turn out to be very instrumental to derive, in the subsequent sections, a number of key conditions for characterizations of both exponential and asymptotic stabilizability of nonlinear control systems.


\section{Exponential Stabilizability of Control Systems via Composition Operators}\label{sec3}


Let us first recall the main result by Gupta, Jafari, Kipka and Mordukhovich \cite[Theorem~3.2]{GJKM}. Given a mapping $f\colon\R^n\times\R^m\to\R^l$ which  is continuously differentiable on a neighborhood $\mcal{U}\times\mcal{V}$ of the origin $(0,0)\in\R^n\times\R^m$, denote the partial Jacobian matrices of $f$ at $(0,0)$ by
\begin{align}
    A_f := \restr{\frac{\partial f}{\partial x}}{(0,0)} \es \textrm{and} \es B_f := \restr{\frac{\partial f}{\partial u}}{(0,0)}. \label{AB def}
\end{align}

\noindent Given further a linear operator/matrix $T\colon\R^n\to\R^n$, we denote its (complex) {\it spectrum} by $\Lambda(T)$ and consider the eigenvalues with {\it nonnegative real parts}

$$\Lambda_+(T) := \cb{\lambda \in \Lambda(T)\;\big|\;Re(\lambda)\geq 0}.$$

\noindent The relevant (for our purposes) portion of \cite[Theorem~3.2]{GJKM} is formulated as follows.

\begin{thm} \label{GJKM 3.2}
    \noindent Let $f$ be of class $C^1$ in a neighborhood $\mcal{U} \times \mcal{V} \subseteq \R^n \times \R^m$ of the origin $(0,0)$, and let the origin be an equilibrium point of $f$. Assume also that $\Lambda_+(A_f) \subseteq \R$, that the mapping $f$ is linearly open around the origin, and that the condition
    \begin{equation}\label{spect}
    cov f(0,0) > \eta_f := {\rm sup}\big\{\lambda\;\big|\;\lambda \in \Lambda_+(A_f)\big\}
    \end{equation}
    is satisfied. Then system \eqref{sys} is locally exponentially stabilizable by means of continuous stationary feedback laws.
    \noindent Furthermore, the linear openness of $f$ around the origin is necessary for such a local exponential stabilization without any additional assumptions on the system data if continuously differentiable stationary feedback laws are used instead.
\end{thm}


The message of Theorem~\ref{GJKM 3.2} can be summarized in the following way: for a vector field $f$ of class $C^1$ in a neighborhood $\mcal{U} \times \mcal{V}$ of $(0,0)\in\R^n\times\R^n$, a certain openness property with respect to the {\it real spectrum} of $A_f$ given by  \eqref{spect} determines the local exponential stabilizability of \eqref{sys} when {\it $C^1$ stationary feedback laws} are required.  The necessity given parallels the known condition of Zabczyk \cite{zab} that \eqref{sys} is locally exponentially stabilizable by $C^1$ feedback laws if and only if its linearization is locally exponentially stabilizable. The sufficiency of this condition parallels the well-known Hautus lemma.  Having this in mind, let us proceed to establish a method by which linear openness may be further applied to the problem of stabilizability. To accomplish this, we develop a novel {\it composition operator approach} that allows us to obtain a partial characterization of exponential stabilizability by means of continuous stationary feedback laws for systems \eqref{sys} which are either linearly open at an equilibrium point or are linearly open and satisfy a diffeo-geometric transversality condition near an equilibrium point.  While these results are certainly implied by Zabczyk's necessary and sufficient condition for exponential stabilizability by stationary $C^1$ feedback laws, they are instrumental in proving the final results concerning stabilizability by merely continuous feedback laws.  \vspace*{0.05in}

To begin with, let us recall some facts on {\it composition operators} of the type
\begin{equation}\label{comp}
    T_h: C^1(\mcal{U} \times \mcal{V},\R^n) \rightarrow C^1(\mcal{O},\R^n)\;\mbox{ with }\;
    T_h f: = f \circ h.
\end{equation}
\noindent
It follows from \eqref{comp} that $h \in C^1(\mcal{O},\mcal{U} \times \mcal{V})$), and we get by the classical chain rule that

\begin{align}
    \restr{J_{T_hf}}{(x,u)} = \restr{J_{f}}{h(x,u)} \restr{J_h}{(x,u)} \label{jac comp}.
\end{align}
\noindent
Writing now $h(x,u)=(h_1(x,u),h_2(x,u))$ with $h_1 \in C^1(\mcal{O}, \mcal{U})$ and $h_2 \in C^1(\mcal{O}, \mcal{V})$, gives us

\begin{align*}
    \restr{J_{T_hf}}{(x,u)} &=  \restr{J_{f}}{h(x,u)} \restr{J_h}{(x,u)} \\
                            &=  \begin{bmatrix}
                            	\restr{\frac{\partial f}{\partial x}}{h(x,u)}\restr{\frac{\partial f}{\partial u}}{h(x,u)}
                            \end{bmatrix}\begin{bmatrix}
                            \restr{J_{h_1}}{(x,u)} \\
                            \restr{J_{h_2}}{(x,u)}
                            \end{bmatrix} \\
                           &= \restr{\frac{\partial f}{\partial x}}{h(x,u)}\restr{J_{h_1}}{(x,u)}+\restr{\frac{\partial f}{\partial u}}{h(x,u)}\restr{J_{h_2}}{(x,u)}.
\end{align*}

This allows us to estimate from both above and below the exact linear openness bound of  composition operators. We say that a mapping $h\colon\R^n\times\R^m\to\R^n\times\R^m$ is {\it stationary} if $h(0,0)=(0,0)$.

\begin{lemma} \label{cov inequal comp op}
    \noindent Let $f$ be of class $C^1$ in a neighborhood $\mcal{U} \times \mcal{V} \subseteq \R^n \times \R^m$ of the origin, let $\mcal{O} \subseteq \R^n \times \R^m$ be a neighborhood of the origin as well, and let $h \in C^1(\mcal{O},\mcal{U} \times \mcal{V})$ be stationary. Then we have the estimates
    \begin{equation}\label{cov-comp}
        \|J_f(0,0)\|\cdot cov (h)(0,0) \geq cov (T_hf)(0,0) \geq cov(f)(0,0)\cdot cov(h)(0,0).
    \end{equation}
\end{lemma}


\begin{proof}
    Under the assumptions made we get from Theorem~\ref{lo fr jac} and the chain rule equalities that
    \begin{align*}
        cov (T_hf)(0,0) &= \textrm{min} \left \{ \left \| \restr{J_{T_hf}^*}{(0,0)} v \right \| : \|v\| = 1 \right \} \\
        &= \textrm{min} \left \{ \left \| \left ( \restr{J_f}{(h(0,0))}\restr{J_h}{(0,0)} \right )^* v \right \| : \|v\| = 1 \right \} \\
        &= \textrm{min} \left \{ \left \| \restr{J_h^*}{(0,0)} \restr{J_f^*}{(h(0,0))} v \right \| : \|v\| = 1 \right \} \\
        &= \textrm{min} \left \{ \left \| \restr{J_h^*}{(0,0)} \restr{J_f^*}{(0,0)} v  \right \| : \|v\| = 1 \right \} \\
        &= \textrm{min} \left \{ \left \| \restr{J_h^*}{(0,0)}  w \right \| : w \in \restr{J_f^*}{(0,0)} \cp{\mbb{S}_{\R^{n+m}}} \right \} \\
        &= \textrm{min} \left \{ \left \| \restr{J_h^*}{(0,0)}  w \right \| : \Gamma \left ( \restr{J_f}{(0,0)} \right ) \leq \| w \| \leq \left \| \restr{J_f}{(0,0)} \right \| \right \}.
    \end{align*}
    
    Therefore, by rescaling $w$ accordingly in the last inequality, we conclude that both estimates in \eqref{cov-comp} hold, which completes the proof of the lemma.
\end{proof}


Combining now Lemmas~\ref{lomr identity} and \ref{cov inequal comp op} leads us to the following interplay between given moduli of linear openness and metric regularity.

\begin{cor} \label{cov suf comp op}
    \noindent Let $f$ be of class $C^1$ in a neighborhood $\mcal{U} \times \mcal{V} \subseteq \R^n \times \R^m$ of the origin, and let $0 < cov(f)(0,0) < \kappa$ for some fixed constant $\kappa < \infty$. Suppose further that $h: \mcal{O} \rightarrow \mcal{U} \times \mcal{V}$ is stationary and of class $C^1$ in a neighborhood $\mcal{O}$ of the origin, and that $cov(h)(0,0) > 0$. Then the estimate $cov(h)(0,0) \geq \kappa reg(f)(0,0)$ implies that $cov(T_hf)(0,0) > \kappa$ for the composition operator \eqref{comp}. Similarly, it follows from $reg (h)(0,0) \kappa \geq \|J_f(0,0)\|$ that $cov(T_hf) \leq \kappa$.
\end{cor}
\begin{proof} 
    Observe that if for some $\kappa < \infty$ we have $cov(f)(0,0) < \kappa$, then $1 < \kappa \ reg(f)(0,0)$ by Lemma~\ref{lomr identity}. Then the claimed assertions follow from Lemma~\ref{cov inequal comp op}.
\end{proof}


\noindent To proceed further with establishing composition extensions of Theorem~\ref{GJKM 3.2}, we need the following two lemmas. The first one is a straightforward computation of the Jacobian matrix, while the second lemma lies at the heart of the main theorem of this section that follows. We say that a mapping is a {\it stationary diffeomorphism} if it is simultaneously stationary and a diffeomorphism.

\begin{lemma}\label{comp jac formula}
    \noindent Suppose that $f$ is of class $C^1$ in a neighborhood $\mcal{U} \times \mcal{V} \subseteq \R^n \times \R^m$ of the origin, which is an equilibrium point of $f$. Then, given a composition operator $T_h:C^1(\mcal{U}\times\mcal{V},\R^n) \rightarrow C^1(\mcal{O},\R^n)$ with $h\in C^1(\mcal{O},\mcal{U} \times \mcal{V})$ being stationary and $C^1$ on a neighborhood $\mcal{O} \subseteq \mcal{U} \times \mcal{V}$ of the origin, we get
    $$\restr{J_{T_hf}}{(0,0)} = \ba{c|c}{
    \cp{\restr{J_f}{(0,0)} A_h} & \cp{\restr{J_f}{(0,0)} B_h}
    },$$
    \noindent where $A$ and $B$ are taken from \eqref{AB def}, and where $\ba{c|c}{X & Y}$ denotes the matrix concatenation of the matrices $X$ and $Y$, i.e., the corresponding block matrix.
\end{lemma}


\begin{proof} 
    Consider the composition operator $T_h f$ from \eqref{comp} under the assumptions made. Then it follows by the direct application of (\ref{jac comp}) that
    \begin{align*}
    \restr{J_{T_hf}}{(0,0)} &=  \restr{J_{f}}{h(0,0)} \restr{J_h}{(0,0)} \\
                            &=  \ba{c|c}{
                            	A_f & B_f
                            }\ba{c|c}{
                            A_{h_1} & B_{h_1} \\
                            \hline
                            A_{h_2} & B_{h_2}
                            } \\
                           &= \ba{c|c}{
                              \cp{A_f A_{h_1}+B_f A_{h_2}} &
                              \cp{A_f B_{h_1}+B_f B_{h_2}}
                           }.
    \end{align*}
    \noindent This can be simply rewritten as
    \begin{align*}
    \ba{c|c}{
    \cp{A_{T_hf}} & \cp{B_{T_hf}}
    } = \ba{c|c}{
    \cp{\restr{J_f}{(0,0)} A_h} & \cp{\restr{J_f}{(0,0)} B_h}
    },
    \end{align*}
    which verifies the claimed Jacobian representation.
\end{proof}


The next lemma allows us to {\em shift positive real eigenvalues} of the spatial portion of the linearization $A_f$ of a linearly open smooth vector field $f$ by  {\em composition} with a {\em local diffeomorphism}.

\begin{lemma}\label{dif comp}
    \noindent Let $f$ be of class $C^1$ in a neighborhood $\mcal{U} \times \mcal{V} \subseteq \R^n \times \R^m$ of the origin which is an equilibrium point of $f$, and let $f$ be linearly open around the origin, i.e., $cov f(0,0)>0$. Then there exists a neighborhood of the origin $\mcal{O} \subseteq \mcal{U}\times\mcal{V}$ and a composition operator $T_h:C^1(\mcal{U}\times\mcal{V},\R^n) \rightarrow C^1(\mcal{O},\R^n)$ such that the mapping $h \in C^1(\mcal{O},\mcal{U} \times \mcal{V})$ in \eqref{comp} is a stationary diffeomorphism with $\Lambda_+\cp{A_{T_hf}}\subseteq \R$.
\end{lemma}


\begin{proof} 
    Under the assumptions made, Theorem~\ref{lo fr jac} tells us that the Jacobian matrix $\restr{J_f}{(0,0)}$ has full rank. Furthermore, it follows from Lemma~\ref{comp jac formula} that
    \begin{align*}
        \ba{c|c}{
        \cp{A_{T_hf}} & \cp{B_{T_hf}}
        } = \ba{c|c}{
        \cp{\restr{J_f}{(0,0)} A_h} & \cp{\restr{J_f}{(0,0)} B_h}
        }.
    \end{align*}
    Denoting now by $A^+$ the right Moore-Penrose pseudoinverse of the linear operator $A$, we get
    \begin{align}
        A_h = -\restr{J_f}{(0,0)}^+ := -\restr{J_f}{(0,0)}^*\left (\restr{J_f}{(0,0)}\restr{J_f}{(0,0)}^* \right)^{-1} \label{pseudo}
    \end{align}
    since $h$ is stationary. It follows from the properties of the Moore-Penrose pseudoinverse that
    
    $$A_{T_hf} = -\restr{J_f}{(0,0)}\restr{J_f}{(0,0)}^+ = -I,$$
    and so $\eta_{T_hf}:=\sup{\lambda\;\big|\;\lambda \in \Lambda_+(A_{T_hf})} = -\infty$
    by the convention that $\sup{\emptyset} = -\infty$. This allows us to conclude that the existence of a stationary mapping $h \in C^1(\mcal{O},\mcal{U} \times \mcal{V})$ in \eqref{comp} satisfying (\ref{pseudo}) yields the validity of all the assertions of the theorem.
    
    Let us now explicitly construct a mapping $h$ with the desired properties. Write $h(x,u) = (h_1(x,u),h_2(x,u))$ and denote by $\cp{Q}_{i,j}$ the $(i,j)^{th}$ entry of a given matrix $Q$. Then define
    
    \begin{align}
        \cp{h_1(x,u)}_{i,1}:&= \sum_{j=1}^{n} \cp{\restr{J_f}{(0,0)}^+}_{i,j}x_j,\quad i=1,\ldots,n, \label{h1 construction}
    \end{align}
    \begin{align}
        \cp{h_2(x,u)}_{k,1}: &= \sum_{j=1}^{n} \cp{\restr{J_f}{(0,0)}^+}_{(n+k),j}x_j + u_k,\quad k=1,\ldots,m, \label{h2 construction}
    \end{align}
    
    More concisely, we set $h(x,u)=\restr{J_f}{(0,0)}^+x+\ba{c}{0\\u}$.  Observe that the operator $h(x,u)$ is linear and stationary. Furthermore, it follows from (\ref{h1 construction}) and (\ref{h2 construction}) that  $A_{h}=\restr{J_f}{(0,0)}^+$ for the partial Jacobian matrix of $h$, and that (\ref{pseudo}) is satisfied. To get a neighborhood $\mcal{O}\subseteq \mcal{U}\times\mcal{V}$ with $h(x,u) \in \mcal{U} \times \mcal{V}$ for all $(x,u) \in \mcal{N}$, we set $\mcal{O}:=r({\rm int}\,\mbb{B}_{\R^n\times \R^m})$, where the radius $r$ is selected as $r: = \sup{r\;\big|\;r\|h\| \leq \sup{t>0\;\big|\;t\mbb{B}_{\R^n\times \R^m} \subseteq\mcal{U}\times\mcal{V}}}$. Since $ker(T^+)=ker(T^*)$ for a linear operator $T$, we deduce that
    \begin{align*}
    ker(A_h) &= ker\left (\restr{J_f}{(0,0)}^+\right ) = ker\left (\restr{J_f^*}{(0,0)}\right ) = \left (rge\left (\restr{J_f}{(0,0)}\right ) \right )^\perp = 0,
    \end{align*}
    where $ker$ and $rge$ denote, respectively, the kernel and range of a linear operator. This clearly implies that $A_h$ is of full rank, and that $h$ is a diffeomorphism due to $h_2(0,u)=u$ as follows from (\ref{h2 construction}). Thus we get the mapping $h$ in composition \eqref{comp} satisfying all the required properties and such that
    $A_{T_hf}$ has no eigenvalues with positive real parts.
\end{proof}

Note that the diffeomorphism constructed in the proof of Lemma~\ref{dif comp} is in fact {\em linear}. While this is convenient for our purposes, nonlinear diffeomorphisms suffice as well provided that they satisfy the equalities in (\ref{pseudo}). \vspace*{0.05in}

Now we are ready to derive a significant extension of the sufficient condition for stabilizability in Theorem~\ref{GJKM 3.2} by showing that the usage of a composition operator with a linear diffeormorphism allows us to {\em completely avoid} the restrictive spectral assumption \eqref{spect}.

\begin{thm}\label{les dif comp} 
    Consider the control system \eqref{sys} in the setting of Lemma~{\rm\ref{dif comp}}. Then there exists a neighborhood of the origin $\mcal{O}\subseteq \mcal{U} \times \mcal{V}$ and a composition operator $T_h:C^1(\mcal{U}\times\mcal{V},\R^n) \rightarrow C^1(\mcal{O},\R^n)$ with a stationary linear diffeomorphism $h\colon{\cal O}\to{\cal U}\times{\cal V}$ such that the modified system
    \begin{equation}\label{modi}
    \dot{x}=T_hf(x,u)\;\mbox{ with }\;T_h f: = f \circ h,\quad t\ge 0,
    \end{equation}
    is locally exponentially stabilizable by means of continuous stationary feedback laws.  Further, it suffices to choose $h$ as in (\ref{h1 construction}) and (\ref{h2 construction}), and the trivial feedback law $u(x)=0$ yields the desired exponential stability for (\ref{modi}).
\end{thm}

\begin{proof} 
    Under the assumptions imposed in this theorem we deduce from Theorem \ref{lo fr jac} that the Jacobian matrix $\restr{J_f}{(0,0)}$ has full rank. Take now the stationary linear operator $T_h$ from \eqref{comp} with the linear diffeomorphism $h$ constructed in the proof of Lemma~\ref{dif comp} and observe similarly to the proof above that $\eta_{T_hf}=\sup{\lambda\;\big|\;\lambda\in \Lambda_+(A_{T_hf})}=-\infty$ and that $\textrm{rank}\cp{A_{T_hf}} = n$ due to $A_{T_hf} = -I$. Thus $\textrm{rank}\cp{\restr{J_{T_hf}}{(0,0)}} \geq n$ by Lemma~\ref{comp jac formula}. Using again Theorem~\ref{lo fr jac} tells us that $cov T_hf(0,0) > 0$ and therefore $cov T_hf(0,0) > \eta_{T_hf}$. Then we get by Theorem~\ref{GJKM 3.2} that the composite system $\dot{x}=T_hf(x,u)$ with the explicitly constructed stationary linear diffeomorphism $h$ is locally exponentially stabilizable by means of continuous stationary feedback laws.  Further, it is routine to check via linearization that the trivial feedback law $u(x)=0$ stabilizes $\dot{x}=T_hf(x,u)$ when $h$ is chosen as in (\ref{h1 construction}) and (\ref{h2 construction}).
\end{proof}


Loosely speaking, Theorem~\ref{les dif comp} can be understood as follows. If the original system \eqref{sys} is not locally exponentially stabilizable by means of stationary continuous feedback laws, then the nice composition operator constructed above provides a "correcting lens" of sorts. which allows $f$ to interpret its domain properly and to achieve a relatively high degree of stability. Let us illustrate this procedure by the following example taken from Sontag in \cite{sontag2} who designed it to show that Brockett's necessary condition from Theorem~\ref{brocketts} is not sufficient for asymptotic stabilizability of \eqref{sys} by using stationary continuously differentiable feedback laws.

\begin{ex}\label{ex 1}  
    Consider system \eqref{sys} with $f(x,u):= x+u^3$. By a routine calculation we can check that the feedback law $u(x):=(-2x)^{1/3}$ locally exponentially stabilizes the system, and it is indeed continuous while not $C^1$ around the origin. At the same time we have
    $$\restr{J_f}{(0,0)} = \ba{c|c}{1 & 0}$$
    which readily implies that the relationships
    $$
    {\rm sup}\big\{\lambda\;\big|\;\lambda \in \Lambda_+(A_f)\big\} =1 \;\mbox{ and }\;cov f(0,0) = 1.
    $$
    We get therefore that $cov f(0,0)= {\rm sup}\{\lambda\;\big|\;\lambda \in \Lambda_+(A_f)\}$, which shows that the spectral sufficient condition \eqref{spect} of Theorem~\ref{GJKM 3.2}
    fails and hence that result does not ensure the local exponential stabilizability of \eqref{sys} by means of stationary, continuously differentiable feedback laws. Moreover, since $f=x+u^3$
    is linearly open at $(0,0)$, the necessity part of Theorem~\ref{GJKM 3.2} tells us that the original system \eqref{sys} is not exponentially stabilizable by means of continuously differentiable stationary feedback laws.
    
    However, the situation dramatically changes when we employ an appropriate composition operator. Observing that
    $$\restr{J_f}{(0,0)}^+ = -\restr{J_f}{(0,0)}^*\left ( \restr{J_f}{(0,0)}\restr{J_f}{(0,0)}^* \right )^{-1} = \begin{bmatrix}
    -1 \\ 0
    \end{bmatrix}$$
    and using $h(x,u): = (-x,u)$ as prescribed by Theorem~\ref{dif comp}, we get $T_hf(x,u)=-x+u^3$ for the composition operator \eqref{comp} ensuring the conditions
    $$\restr{J_{T_hf}}{(0,0)} =\ba{c|c}{-1 & 0},\;cov f(0,0)=1,\;\mbox{ and }\;{\rm sup}\big\{\lambda\;\big|\;\lambda \in \Lambda_+(A_{T_hf})\big\}= -\infty.
    $$
    Employing Theorem~\ref{les dif comp} allows us to conclude that the modified system \eqref{modi} is locally exponentially stabilizable by $C^1$ feedback laws. Indeed, the feedback control achieving the exponential stabilization is the trivial one $\bar u(x)=0$.
\end{ex}


\section{Transversality and Characterizations of Stabilizability}\label{sec4}


In this section we continue our study of feedback stabilizability for nonlinear control systems. Theorem~\ref{les dif comp} suggests to consider, along with the original control system \eqref{sys}, its modified version \eqref{modi}, which allows us to avoid the restrictive spectral assumption imposed in Theorem~\ref{GJKM 3.2}. Now we proceed with further extensions, which eventually lead us to some further characterizations of both local exponential and asymptotic feedback stabilizability of linearly open systems. First, however, let us consider a motivating example.

\begin{ex}\label{no cov no exp stab}  
    Consider the system \eqref{sys}  constructed by Coron in \cite{coron2}, where $f$ is given by:
    
    $$f(x,u): = \begin{bmatrix}
    x_2^3 - 3(x_1-x_3)^2x_2 \\
    (x_1-x_3)^3-3(x_1-x_3)^2x_2\\
    u
    \end{bmatrix}.$$
    As shown in \cite{coron2}, while this system satisfies Brockett's necessary condition, it is not locally asymptotically stabilizable. However, this example does not contradict our preceding results since the above system is {\em not linearly open} around equilibrium points. Indeed, observe first that all the equilibria of $f$ are given by $x_1 = x_3$, $x_2=0$, and $u=0$. It is straightforward to check that at any equilibrium point $x_e$ we have
    
    $$\restr{J_f}{(x_e)} = \ba{ccc|c}{0 & 0 & 0 & 0 \\ 0 & 0 & 0 & 0 \\ 0 & 0 & 0 & 1}.$$
    
    \noindent Thus $cov f(x_e) = 0$ by Theorem~\ref{lo fr jac}, which confirms that $f$ fails to be linearly open around $x_e$.  Similar examples on Brockett's theorem can be found in \cite{coron1,sontag3}. Routine calculations allow us to check that in all these examples vector fields are not linearly open around equilibrium points.\vspace*{0.05in}
    However, after some straightforward computation, we have $\restr{\frac{\partial f}{\partial u}}{(x,u)} = \ba{ccc}{0 & 0 & 1}^*$ and
    
    $$\restr{\frac{\partial f}{\partial x}}{(x,u)} = \ba{ccc}{-6x_2(x_1-x_3) & 3x_2^2-3(x_1-x_3)^2 & 6x_2(x_1-x_3) \\ 3(x_1-x_3)^2-6x_2(x_1-x_3) & -3(x_1-x_3)^2 & 6x_2(x_1-x_3)-3(x_1-x_3)^2 \\ 0 & 0 & 0}.$$
    
    \noindent This demonstrates that the tangent space satisfies the transversality condition
    $\frac{\partial f}{\partial x_1} =-\frac{\partial f}{\partial x_3}$ in the $(x_1,x_3)$-plane at the equilibrium point $x_1=x_3,\;x_2=0,\;u=0$. However, the term $\frac{\partial f}{\partial x_2} $ approaches the equilibrium tangentially at a higher order, and in the other coordinates it fails to satisfy any transversality condition at all. Such a local analysis near equilibrium suggests that a linear approach to an equilibrium, coupled with a certain transversality condition is needed to guarantee exponential stabilizability (indeed, similar observations were made in, e.g., \cite{andreini}).

\end{ex}


Along this line, it is worth noting that Brockett's necessary condition is a \textit{projectability} condition yielding (at least locally) a well-defined pushforward. On the other hand, an appropriate transversality constraint provides a way to understand how the tangent manifolds approach the equilibrium.  In what follows, we will show that this, in fact, is precisely the case.  Let us now define several concepts which will be of central importance for the subsequent material. We use them heavily and formalize them as follows.

\begin{defi}\label{Tdef}\

    \noindent Suppose $f$ is of class $C^1$ on some open set $\mcal{D}\subseteq \R^n\times\R^m$, and let $(a,b)\in \mcal{D}$.
    
    \begin{enumerate}[(i)]
        \item If $ n = m $, we say $f$ \textit{satisfies a transversality condition at} $(a,b)$ if there exists $c \geq 0$ such that $\restr{\frac{\partial f}{\partial x}}{(a,b)} = -c\restr{\frac{\partial f}{\partial u}}{(a,b)}$. Further, we say that $c$ is the \textit{transversality constant of} $f$ \textit{at} $(a,b)$. \label{TC}
        \item We say $f$ \textit{satisfies a type I semitransversality condition at} $(a,b)$ if there exists an $m\times n$ matrix $Q$ such that $\restr{\frac{\partial f}{\partial x}}{(a,b)} = -\restr{\frac{\partial f}{\partial u}}{(a,b)}Q$.  Further, we call any matrix $Q$ satisfying the above a \textit{transverse factor} of $f$ \textit{at} $(a,b)$. \label{STC}
        \item We say $f$ \textit{satisfies a type II semitransversality condition at} $(a,b)$ if $$\mbb{P}_{ker\left(\restr{\frac{\partial f}{\partial u}}{(a,b)}^*\right)^\perp}\left(\restr{\frac{\partial f}{\partial x}}{(a,b)}\right)\left(\restr{\frac{\partial f}{\partial x}}{(a,b)}\right)^*\mbb{P}_{ker\left(\restr{\frac{\partial f}{\partial u}}{(a,b)}^*\right)^\perp}=\left(\restr{\frac{\partial f}{\partial x}}{(a,b)}\right)\left(\restr{\frac{\partial f}{\partial x}}{(a,b)}\right)^*.$$ \label{QTC}
    \end{enumerate}

\end{defi}


Clearly, (\ref{TC}) implies (\ref{STC}) and--noting that $X=M^+N$ is a solution to the matrix equation $MX=N$ if any such solutions exist--it follows from $MM^+=\mbb{P}_{ker(M^*)^\perp}$ that (\ref{STC}) implies (\ref{QTC}).  The next theorem serves as a simple introduction to the manner in which composition operators may be used to demonstrate stability for linearly open systems satisfying some form of a transversality condition.  However, before beginning, it is important to note that, while portions of several of the preliminary results in what follows may be derived from the well-known Hautus lemma, the final results obtained at the end of the section certainly exceed the Hautus lemma's scope.

\begin{thm} \label{T,I}

    \noindent Suppose $f$ is of class $C^1$ in a neighborhood $\mcal{U} \times \mcal{V} \subseteq \R^n \times \R^n$ of the origin, and let the origin be an equilibrium point of $f$.  Suppose $cov f(0,0)>0$, and suppose $f$ satisfies a transversality condition at the origin.  Then \eqref{sys} is locally exponentially stabilizable by means of continuously differentiable, stationary feedback laws.  Moreover, letting $c$ denote the transversality constant of $f$ at the origin, this stability is achieved by the linear feedback law
    
    $$u(x)=-(c^2+1)B_f^*\cp{A_fA_f^*+B_fB_f^*}^{-1}(I+A_f)x.$$

\end{thm}


\begin{proof}
    Suppose $f$ is of class $C^1$ in a neighborhood $\mcal{U} \times \mcal{V} \subseteq \R^n \times \R^n$ of the origin, and let the origin be an equilibrium point of $f$.  Suppose $cov f(0,0)>0$, and suppose $f$ satisfies a transversality condition at the origin.  That is, for some $c\geq0$, suppose $A_f = -cB_f$.  By the assertion that $covf(0,0)>0$, we may conclude that $\restr{J_f}{(0,0)}$ has full row-rank.  Correspondingly, $\restr{J_f}{(0,0)}$ has a right Moore-Penrose pseudoinverse.  As such, let $T_h$ be the composition operator given in (\ref{h1 construction}) and (\ref{h2 construction}).  Let $F\in C^1(\R^n\times\R^n,\R^n)$ be given by $F(x,y) = T_hf\cp{x,0}-f\cp{x,(c^2+1)B_f^*y}$.  Observe that $F(0,0)=0$, and
    
    $$\restr{\frac{\partial F}{\partial x}}{(0,0)}=-I-A_f, \quad \restr{\frac{\partial F}{\partial y}}{(0,0)}=-(c^2+1)B_fB_f^*=-(A_fA_f^*+B_fB_f^*).$$
    
    As $covf(0,0)>0$ and as $f$ satisfies a transversality condition at the origin, it follows that either $c=0$ and $rank(B_f)=n$, or $c \neq 0$ and $rank(A_f)=rank(B_f)=n$.  In either case, $A_fA_f^*+B_fB_f^*$ is invertible, so by the Implicit Function Theorem, there exists an open set $\mcal{U}_1 \subseteq \mcal{U}$ containing the origin and a unique function $g \in C^1(\mcal{U}_1,\mcal{U})$ such that $g(0)=0$ and, for all $x \in \mcal{U}_1$, $F\cp{x,g(x)} = 0$.  That is to say, on $\mcal{U}_1$, we have $T_hf\cp{x,0}=f\cp{x,(c^2+1)B_f^*g(x)}$.
    
    So, by the exponential stability of $\dot{x}=T_hf(x,0)$, the system, $\dot{x}=f\cp{x,(c^2+1)B_f^*g(x)}$ is locally exponentially stable and the system trajectories $\phi(t,x_0)\rightarrow 0$ for all $x_0 \in \mcal{U}_1$. Moreover,  by the implicit function theorem, $g$ satisfies
    
    \begin{align*}
        \restr{\frac{\partial g}{\partial x}}{0} &= -\left (\restr{\frac{\partial F}{\partial y}}{(0,g(0))} \right )^{-1}\restr{\frac{\partial F}{\partial x}}{(0,0)} = -\cp{A_fA_f^*+B_fB_f^*}^{-1}(I+A_f).
    \end{align*}
    
    As such, it is routine to check that $u(x)=(c^2+1)B_f^*\restr{\frac{\partial g}{\partial x}}{0}x$ is a continuously differentiable stationary control stabilizing the system.
\end{proof}


While the transversality condition (\ref{TC}) of Definition~\ref{Tdef} used in Theorem~\ref{T,I} is relatively strong, we will now show that the much weaker type I semitransversality condition (\ref{STC}) of Definition~\ref{Tdef} is indeed sufficient for the same conclusion.  Note that the result above, as well as those below, concern exponential feedback stabilizability of the {\it original linearly open} control system \eqref{sys} without imposing the spectral assumption \eqref{spect}.

\begin{thm}\label{ST,I}

    \noindent Suppose $f$ is of class $C^1$ in a neighborhood $\mcal{U} \times \mcal{V} \subseteq \R^n \times \R^m$ of the origin, and let the origin be an equilibrium point of $f$.  Suppose $cov f(0,0)>0$, and suppose $f$ satisfies a type I semitransversality condition at the origin.  Then \eqref{sys} is locally exponentially stabilizable by means of continuously differentiable, stationary feedback laws.  Moreover, if $Q$ is any transverse factor of $f$ at the origin, this stability is achieved by the linear feedback law
    
    $$u(x)=-(QQ^*+I)B_f^*\left(A_fA_f^*+B_fB_f^*\right)^{-1}\left(I+A_f\right)x.$$

\end{thm}


\begin{proof}
    Suppose $f$ is of class $C^1$ in a neighborhood $\mcal{U} \times \mcal{V} \subseteq \R^n \times \R^m$ of the origin, and let the origin be an equilibrium point of $f$.  Suppose $cov f(0,0)>0$, and suppose $f$ satisfies the type I semitransversality condition $A_f = -B_fQ$ for some $m\times n$ matrix $Q$. By the assertion that $covf(0,0)>0$, we may conclude that $\restr{J_f}{(0,0)}$ has full row-rank.  Correspondingly, $\restr{J_f}{(0,0)}$ has a right Moore-Penrose pseudoinverse.  As such, let $T_h$ be the composition operator given in (\ref{h1 construction}) and (\ref{h2 construction}).  Let $F\in C^1(\R^n\times\R^n,\R^n)$ be given by
    
    $$F(x,y) = T_hf\cp{x,0}-f\cp{x,(QQ^*+I)B_f^*y}.$$
    
    Observe that $F(0,0)=0$, and note that $\restr{\frac{\partial F}{\partial x}}{(0,0)} = -\left (I+A_f \right )$.  Similarly,
    
    $$\restr{\frac{\partial F}{\partial y}}{(0,0)} = -B_f(QQ^*+I)B_f^* = - \left( B_fQQ^*B_f^*+B_fB_f^*\right) = -\left(A_fA_f^*+B_fB_f^* \right ).$$
    
    As $covf(0,0)>0$, it follows that $\restr{J_f}{(0,0)}$ is full rank and $\restr{J_f}{(0,0)}\restr{J_f}{(0,0)}^*=A_fA_f^*+B_fB_f^*$ is invertible.  So, by the Implicit Function Theorem, there exists an open set $\mcal{U}_1 \subseteq \mcal{U}$ containing the origin and a unique function $g \in C^1(\mcal{U}_1,\mcal{U})$ such that $g(0)=0$ and, for all $x \in \mcal{U}_1$, $F\cp{x,g(x)} = 0$.  That is to say, on $\mcal{U}_1$, we have
    
    $$T_hf\cp{x,0}=f\cp{x,(QQ^*+I)B_f^*g(x)}.$$
    
    So, by the stability of $\dot{x}=T_hf(x,0)$, the system $\dot{x}=f\cp{x,(QQ^*+I)B_f^*g(x)}$ is locally exponentially stable and $\phi(t,x_0)\rightarrow 0$ for all $x_0 \in \mcal{U}_1$.  Moreover, by the implicit function theorem, $g$ satisfies
    
    \begin{align*}
        \restr{\frac{\partial g}{\partial x}}{0} &= -\left (\restr{\frac{\partial F}{\partial y}}{(0,g(0))} \right )^{-1}\restr{\frac{\partial F}{\partial x}}{(0,0)} = -\left(A_fA_f^*+B_fB_f^*\right)^{-1}\left(I+A_f\right).
    \end{align*}
    
    As such, it is routine to check by linearization that $u(x)=-(QQ^*+I)B_f^*\left(A_fA_f^*+B_fB_f^*\right)^{-1}\left(I+A_f\right)x$ is a continuously differentiable stationary control stabilizing the system.

\end{proof}


Following this, we will now show that the type I semitransversality condition (\ref{STC}) of Definiton~\ref{Tdef} used in Theorem~\ref{ST,I} may also be relaxed to the type II transversality condition (\ref{QTC}) of Definiton~\ref{Tdef}.

\begin{thm} \label{QT,I}

    \noindent Suppose $f$ is of class $C^1$ in a neighborhood $\mcal{U} \times \mcal{V} \subseteq \R^n \times \R^m$ of the origin, and let the origin be an equilibrium point of $f$.  Suppose $cov f(0,0)>0$, and suppose $f$ satisfies a quasitransversality condition at the origin.  Then, $\dot{x}=f(x,u)$ is locally exponentially stabilizable by means of continuously differentiable, stationary feedback laws.  Moreover, this stability is achieved by the linear feedback law
    
    $$u(x)=-(B_f^+A_fA_f^*(B_f^+)^*+I)B_f^*\left(A_fA_f^*+B_fB_f^*\right)^{-1}\left(I+A_f\right)x.$$

\end{thm}


\begin{proof}
    Suppose $f$ is of class $C^1$ in a neighborhood $\mcal{U} \times \mcal{V} \subseteq \R^n \times \R^m$ of the origin, and let the origin be an equilibrium point of $f$.  Suppose $cov f(0,0)>0$, and suppose $f$ satisfies the type II semitranversality condition $$\mbb{P}_{ker\left(B_f^*\right)^\perp}A_fA_f^*\mbb{P}_{ker\left(B_f^*\right)^\perp}=A_fA_f^*.$$  By the assertion that $covf(0,0)>0$, we may conclude that $\restr{J_f}{(0,0)}$ has full row-rank.  Correspondingly, $\restr{J_f}{(0,0)}\restr{J_f}{(0,0)}^*$ is invertible.  Noting that $B_fB_f^+ = \mbb{P}_{ker(B_f^*)^\perp}$, observe that
    
    \begin{align*}
        B_f\left (B_f^+A_fA_f^*(B_f^+)^*+I \right )B_f^* &= B_fB_f^+A_fA_f^*(B_f^+)^*B_f^*+B_fB_f^*\\
        &=\mbb{P}_{ker(B_f^*)^\perp}A_fA_f^*\mbb{P}_{ker(B_f^*)^\perp}^*+B_fB_f^* \\
        &= \mbb{P}_{ker(B_f^*)^\perp}A_fA_f^*\mbb{P}_{ker(B_f^*)^\perp}+B_fB_f^*\\
        &=A_fA_f^*+B_fB_f^*
    \end{align*}
    
    As such, linearizing $\dot{x}=f\cp{x,u(x)}$ for the control $u(x)$ given in the statement of the theorem, it is routine to verify that we produce
    
    $$A_f-\left(\restr{J_f}{(0,0)}\restr{J_f}{(0,0)}^*\right)\left(\restr{J_f}{(0,0)}\restr{J_f}{(0,0)}^*\right)^{-1}(I+A_f) = A_f-A_f-I = -I.$$
    
    Hence, the system is locally exponentially stabilizable by continuous, stationary feedback laws, and the control $u(x)$ given in the statement of the theorem is a continuously differentiable stationary control stabilizing the system.

\end{proof}


Following Theorems~\ref{T,I},\ref{ST,I}, and \ref{QT,I}, the sufficiency of both type I and type II semitransversality for local exponential stabilizability of linearly open systems suggests that these properties may warrant further investigation.  The next lemma reveals a necessary and sufficient characterization of type II semitransversality for linearly open systems, and will be instrumental in the full characterization of these related properties.

\begin{lemma}\label{QT,LI}
    Suppose $f$ is of class $C^1$ in a neighborhood $\mcal{U} \times \mcal{V} \subseteq \R^n \times \R^m$ of the origin, and suppose $cov f(0,0)>0$.  Then, $f$ satisfies a type II semitransversality condition at the origin if and only if
    
    $$A_fA_f^* = \frac{1}{2}\left(\mbb{P}_{ker(B_f^*)^\perp}A_fA_f^*+A_fA_f^*\mbb{P}_{ker(B_f^*)^\perp}\right).$$

\end{lemma}


\begin{proof}
    Suppose $f$ is of class $C^1$ in a neighborhood $\mcal{U} \times \mcal{V} \subseteq \R^n \times \R^m$ of the origin, and suppose $cov f(0,0)>0$.  First, observe that if $f$ satisfies a type II semitransversality condition at the origin, then
    
    \begin{align*}
        0&=\left(I - \mbb{P}_{ker(B_f^*)^\perp} \right )\mbb{P}_{ker(B_f^*)^\perp}A_fA_f^*\mbb{P}_{ker(B_f^*)^\perp}\left(I - \mbb{P}_{ker(B_f^*)^\perp} \right ) \\
        &= \left(I - \mbb{P}_{ker(B_f^*)^\perp} \right )A_fA_f^*\left(I - \mbb{P}_{ker(B_f^*)^\perp} \right ) \\
        &= A_fA_f^*-\mbb{P}_{ker(B_f^*)^\perp}A_fA_f^*-A_fA_f^*\mbb{P}_{ker(B_f^*)^\perp} +\mbb{P}_{ker(B_f^*)^\perp}A_fA_f^*\mbb{P}_{ker(B_f^*)^\perp}\\
        &= 2A_fA_f^* - \left(PA_fA_f^*+A_fA_f^*P\right)
    \end{align*}
    
    So, $A_fA_f^* = \frac{1}{2}\left(\mbb{P}_{ker(B_f^*)^\perp}A_fA_f^*+A_fA_f^*\mbb{P}_{ker(B_f^*)^\perp}\right)$ must hold if $f$ satisfies a type II semitransversality condition at the origin.  In the reverse direction, suppose $A_fA_f^* = \frac{1}{2}\left(PA_fA_f^*+A_fA_f^*P\right)$ holds.  Then, by the assumption that $cov f(0,0)>0$, observe that $X=\restr{J_f}{(0,0)}\restr{J_f}{(0,0)}^*$ is a positive definite, symmetric matrix.  Moreover, it is routine to verify that $X$ is a solution to the following discrete time Lyapunov equations for all $c\neq 0$
    
    $$c^2\left(I-\mbb{P}_{ker(B_f^*)^\perp}\right)X\left(I-\mbb{P}_{ker(B_f^*)^\perp}\right)-X+\frac{1}{2}\left(\left(\mbb{P}_{ker(B_f^*)^\perp}\restr{J_f}{(0,0)}\restr{J_f}{(0,0)}^*\right)+\left(\mbb{P}_{ker(B_f^*)^\perp}\restr{J_f}{(0,0)}\restr{J_f}{(0,0)}^*\right)^*\right)=0.$$
    
    Hence, it follows that the system $x_{k+1}=c\left(I-\mbb{P}_{ker(B_f^*)^\perp}\right)x_k$ is globally asymptotically stable for all $c \neq 0$.  But, as a discrete-time linear system of the form $x_{k+1}=Mx_k$ is globally asymptotically stable if and only if $\Lambda(M)\subseteq \{ z \in \mathbb{C}: |z| < 1 \}$, it then must hold that $c(I-\mbb{P}_{ker(B_f^*)^\perp})=0$ (apply, for example, $c=2$).  Hence, $\mbb{P}_{ker(B_f^*)^\perp}=I$ and $\mbb{P}_{ker(B_f^*)^\perp}A_fA_f^*\mbb{P}_{ker(B_f^*)^\perp}=IA_fA_f^*I=A_fA_f^*$. Thus, $f$ satisfies a type II semitransversality condition at the origin if $A_fA_f^* = \frac{1}{2}\left(PA_fA_f^*+A_fA_f^*P\right)$ holds.
\end{proof}


Following Lemma~\ref{QT,LI}, an immediate corollary arises from the global asymptotic stability of the discrete-time dynamical system constructed in the proof.  This reveals that the usually distinct properties of type I and type II semitransversality coincide in the presence of linear openness and are characterized by the rather restrictive row-rank condition imposed on $B_f$.

\begin{lemma}\label{QT,LII}
    Suppose $f$ is of class $C^1$ in a neighborhood $\mcal{U} \times \mcal{V} \subseteq \R^n \times \R^m$ of the origin, and suppose $cov f(0,0)>0$.  Then the following are equivalent:
    
    \begin{enumerate}
        \item $f$ satisfies a type I semitransversality condition at the origin
        \item $f$ satisfies a type II semitransversality condition at the origin
        \item $B_f$ has row-rank $n$.
    \end{enumerate}

\end{lemma}


\begin{proof}
    Suppose $f$ is of class $C^1$ in a neighborhood $\mcal{U} \times \mcal{V} \subseteq \R^n \times \R^m$ of the origin, and suppose $cov f(0,0)>0$.  Suppose $f$ satisfies a type II semitransversality condition at the origin.  Then, $A_fA_f^* = \frac{1}{2}\left(\mbb{P}_{ker(B_f^*)^\perp}A_fA_f^*+A_fA_f^*\mbb{P}_{ker(B_f^*)^\perp}\right)$ holds by Lemma~\ref{QT,LI}.  Moreover, as in the proof of Lemma~\ref{QT,LI}, $\mbb{P}_{ker(B_f^*)^\perp}=I$.  So, as $ker(B_f^*)^\perp=rge(B_f)$, it then follows  that $B_f$ has row-rank $n$.  Hence, $B_f^+=B_f\left(B_fB_f^*\right)^{-1}$, and $B_fB_f^+=I$.  So, clearly, $Q=B_f^+A_f$ is a transverse factor of $f$ at the origin, and $f$ satisfies a type I semitransversality condition at the origin.
\end{proof}


Thus, the type I and type II semitransversality conditions (\ref{STC}) and (\ref{QTC}) of Definition~\ref{Tdef} are \textit{equivalent} in the presence of linear openness.  As such, we will say a linearly open system \textit{satisfies a semitransversality condition} when it satisfies either a type I or type II semitransversality condition.

While the conditions (\ref{TC}), (\ref{STC}) and (\ref{QTC}) of Definiton~\ref{Tdef} at the origin are sufficient for local exponential stabilizability for systems which are linearly open at the origin, straightforward examples such as $f(x,u)=-x$ clearly show that it is not necessary.  However, it is possible to construct a new system--closely related to the original system--that clearly displays the connection between linear openness, the conditions of Definition~\ref{Tdef}, and the stabilizing composition operator constructed in Section~\ref{sec3}.  In the following theorem, we show that any \textit{linearly open} system generates a family of systems parameterized by $C^1(\R^n,\R^n)$, \textit{all} of which are \textit{exponentially stabilizable} and \textit{all} of which \textit{satisfy a transversality condition} (\ref{TC}) at the origin.  Moreover, we show that the composition operator given in Theorem \ref{les dif comp} coincides with the stabilizing control of any system in this family whose parameter in $C^1(\R^n,\R^n)$ is a constant function.

\begin{thm}\label{QTN,I}

    \noindent Suppose $f$ is of class $C^1$ in a neighborhood $\mcal{U} \times \mcal{V} \subseteq \R^n \times \R^m$ of the origin, and let the origin be an equilibrium point of $f$.  Suppose $covf(0,0)>0$.  Let $y \in \R^n$, let $w=\ba{c}{y \\ u}$, and let $G\in C^1\cp{\R^n,\R^n}$.  For $c \neq 0$, define $F:\R^n\times\cp{\R^n\times \R^m} \rightarrow \R^n$ by $F(x,w):=G(x)-G(0)+cf(w)$.  Then, $F$ satisfies a transversality condition at the origin, and the system $\dot{x}=F(x,w)$ is exponentially stabilizable by means of continuously differentiable stationary feedback laws.  Moreover, if $G$ is constant and $c=1$, then the control $w(x)$ stabilizing the system, as given in Theorem~\ref{T,I}, Theorem~\ref{ST,I}, or Theorem~\ref{QT,I} satisfies
    
    $$\cp{x,w(x)} = \cp{x,h(x,0)} = H(x,0)$$
    
    \noindent where $h$ and $H$ are as in (\ref{h1 construction}) and (\ref{h2 construction}) for $f(x,u)$ and $F(x,w)$, respectively.

\end{thm}


\begin{proof}
    Suppose $f$ is of class $C^1$ in a neighborhood $\mcal{U} \times \mcal{V} \subseteq \R^n \times \R^m$ of the origin, and let the origin be an equilibrium point of $f$.  Suppose $covf(0,0)>0$.   Let $y \in \R^n$, let $w=\ba{c}{y \\ u}$, and let $G\in C^1\cp{\R^n,\R^n}$.  For $c \neq 0$, define $F:\R^n\times\cp{\R^n\times \R^m} \rightarrow \R^n$ by $F(x,w):=G(x)-G(0) +cf(w)$.  Then, clearly $\dot{x}=F(x,w)$ has an equilibrium at the origin, and
    
    $$A_F = \restr{\frac{\partial G}{\partial x}}{0}, \es B_F = c\restr{J_{f}}{(0,0)}.$$
    
    Correspondingly, $B_F$ has row-rank $n$ and $covF(0,0)>0$.  So, by Lemma~\ref{QT,LI}, $F$ satisfies a semitransversality condition at the origin.  The local exponential stabilizability of $\dot{x}=F(x,w)$ then follows from Theorem~\ref{QT,I}.  Moreover, if $G$ is constant and $c=1$ and $h$ is as in (\ref{h1 construction}) and (\ref{h2 construction}) for $f$, it follows from Theorem~\ref{les dif comp} that $w(x)=\cp{x,h(x,0)}$ is a stabilizing feedback law for this system, and direct computation verifies that $\cp{x,w(x)} = H(x,0) = \cp{x,h(x,0)}$, where $H$ is as in (\ref{h1 construction}) and (\ref{h2 construction}) for $F$ and $w(x)$ is as in Theorem~\ref{T,I}, Theorem~\ref{ST,I}, or Theorem~\ref{QT,I}.
\end{proof}


So, the existence of a stabilizing composition operator for \eqref{sys} as generated in Section~\ref{sec3} is equivalent to the stabilizability of a modified class of systems induced by the vector field generating the dynamics of \eqref{sys}.  Theorem~\ref{QTN,I} provides additional lines of inquiry that will be discussed in Section~\ref{sec5}.

Next, we will show that the presence of linear openness and a semitransversality condition on domain nearby the origin is sufficient for local exponential stability by continuous feedback laws, allowing us to assess the stability of certain systems with both an equilibrium and a critical point at the origin.  But first, we need a lemma.

\begin{lemma}\label{ST,LIII}
    
    \noindent Suppose $f$ is of class $C^1$ in a neighborhood $\mcal{U} \times \mcal{V} \subseteq \R^n \times \R^m$ of the origin, and let the origin be an equilibrium point of $f$.  Suppose $(x_1,u_1)\in \mcal{U}\times\mcal{V}$, $cov f(x_1,u_1) > 0$, and suppose $f$ satisfies a semitransversality condition at $(x_1,u_1)$.  Define
    
    $$F(x,u):=f(x+x_1,u+u_1)-f(x_1,u_1).$$
    
    \noindent Then, the system $\dot{x}=F(x,u)$ is locally exponentially stabilizable by means of continuously differentiable, stationary feedback laws. Moreover, this stability is achieved by the linear feedback law
    
    $$u(x)=-(B_F^+A_FA_F^*(B_F^+)^*+I)B_F^*\left(A_FA_F^*+B_FB_F^*\right)^{+}\left(I+A_F\right)x.$$

\end{lemma}


\begin{proof}
    Suppose $f$ is of class $C^1$ in a neighborhood $\mcal{U} \times \mcal{V} \subseteq \R^n \times \R^m$ of the origin, and let the origin be an equilibrium point of $f$.  Suppose $cov f(x_1,u_1)>0$, and suppose $f$ satisfies a semitransversality condition condition at $(x_1,u_1)$.
    
    By the assertion that $covf(x_1,u_1)>0$, we may conclude that $\restr{J_f}{(x_1,u_1)}$ has full row-rank.  Define $F(x,u)=f(x+x_1,u+u_1)-f(x_1,u_1)$.  Then, clearly $F(x,u)$ has the origin as an equilibrium, and $\restr{J_F}{(0,0)}=\restr{J_f}{(x_1,u_1)}$.  So, $covF(0,0)>0$ and
    
    $$A_F\mbb{P}_{ker\left(B_F^*\right)^\perp}=\restr{\frac{\partial f}{\partial x}}{(x_1,u_1)}\mbb{P}_{ker\left(\restr{\frac{\partial f}{\partial u}}{(x_1,u_1)}^*\right)^\perp}=\mbb{P}_{ker\left(\restr{\frac{\partial f}{\partial u}}{(x_1,u_1)}^*\right)^\perp}\restr{\frac{\partial f}{\partial u}}{(x_1,u_1)} = \mbb{P}_{ker\left(B_F^*\right)^\perp}B_F.$$
    
    Thus, $F$ satisfies a semitransversality condition at the origin and by Theorem \ref{QT,I}, it follows that the system $\dot{x}=F(x,u)$ is locally exponentially stabilizable by means of continuous stationary feedback laws and that
    
    $$u(x)=-(B_F^+A_FA_F^*(B_F^+)^*+I)B_F^*\left(A_FA_F^*+B_FB_F^*\right)^{+}\left(I+A_F\right)x$$
    
    \noindent is correspondingly a continuously differentiable, stationary stabilizing control.
\end{proof}

Following this, we may produce our aforementioned result, which we first give in the following specific form. The presented characterizations  of the local exponential stabilizability goes beyond standard criteria by involving higher-order terms, while it requires the semitransversality conditions, which are equivalent to the restrictive row-rank property of  $B_f$ due to Lemma~\ref{QT,LII} .

\begin{thm} \label{ST,II}

\noindent Suppose $f$ is of class $C^1$ in a neighborhood $\mcal{U} \times \mcal{V} \subseteq \R^n \times \R^m$ of the origin, and let the origin be an equilibrium point of $f$.  Suppose $\mcal{O}\subseteq \mcal{U}\times\mcal{V}$ is a neighborhood of the origin such that $cov f(x,u) > 0$ for all $(x,u)\in\mcal{O}\setminus\left\{0\right\}$, and suppose  $f$ satisfies a semitransversality condition on $\mcal{O}$.  Then, $\dot{x}=f(x,u)$ is locally exponentially stabilizable by means of continuous, stationary feedback laws which are differentiable on a punctured neighborhood of the origin.

\end{thm}


The proof we will give can be summarized in the following three steps:  First, we use the semitransversality condition to exploit a family of locally exponentially stabilizable systems (as given by Lemma~\ref{ST,LIII}) in such a manner as to produce a criteria under which a Lyapunov function for the system exists if a certain stationary feedback law can be found.  Then, we demonstrate that the existence of such a feedback law is equivalent to a fixed point for a certain nonlinear mapping.  Finally, we show that this mapping does indeed have a fixed point, yielding the desired result.  

To avoid excessive length, we present this argument via a series of relatively self-contained lemmas.  For convenience, let us also adopt the following notation for the remainder of this section: 

For $(y,v) \in \mcal{O}\setminus \left \{ (0,0) \right \}$, write $f[y,v](x,u) = f(x+y,u+v)-f(y,v)$.  For $f$ satisfying the conditions of Lemma~\ref{ST,LIII}) at $(x_1,u_1)=(y,v)$, the systems $\dot{x}=f[y,v](x,u)$ are locally exponentially stabilizable by means of linear feedback laws.  So, in such cases, let us write  $w[y,v](x)$ to denote this linear feedback law.  Similarly, write $F[y,v](x) = f[y,v]\cp{x,w[y,v](x)}$.

\begin{lemma}\label{supporting lemma 1}
    Suppose $f$ satisfies the conditions of Theorem~\ref{ST,II}. Let $\mcal{N} \subseteq \mcal{U}\cap \mcal{O}$ be a compact neighborhood of the origin with nonempty interior.  Then, if there exists a stationary control $u(x)$ which is continuous on $\mcal{N}$ and differentiable on $\mcal{N}\setminus \left \{ 0 \right \}$ such that for all $y \in \mcal{N} \setminus \left \{ 0 \right \}$

    \vspace*{-0.15cm}\begin{align}\restr{\frac{\partial u}{\partial x}}{y} = \left(\restr{\frac{\partial f}{\partial u}}{(y,u(y))}\right)^+\restr{\frac{\partial f}{\partial u}}{(y,u(y))}\restr{\frac{\partial w[y,u(y)]}{\partial x}}{0}, \label{maineq}\end{align} 

    \noindent the system $\dot{x}=f(x,u)$ is locally exponentially stabilizable by means of continuous, stationary feedback laws which are differentiable on a punctured neighborhood of the origin.  
\end{lemma}
\begin{proof}
    Observe that the linearity of $w[y,v]$ and semitransversality of $f$ on $\mcal{O}$ imply
    
    \vspace*{-0.5cm}\begin{align*} 
    \restr{J_{F[y,v]}}{0} &= \restr{\frac{\partial f}{\partial x}}{(y,v)}+\restr{\frac{\partial f}{\partial u}}{(y,v)}\restr{\frac{\partial w[y,v]}{\partial x}}{0} = \restr{\frac{\partial f}{\partial u}}{(y,v)}\left ( \restr{\frac{\partial w[y,v]}{\partial x}}{0} - Q\cp{y,v} \right ).
    \end{align*}
    
    Moreover, as $\restr{J_{F[y,v]}}{0} = -I$, it then follows that 
    
    \vspace*{-0.25cm}$$I=\restr{\frac{\partial f}{\partial u}}{(y,v)}\left ( Q\cp{y,v} - \restr{\frac{\partial w[y,v]}{\partial x}}{0} \right ).$$
    
    Suppose $\mcal{N} \subseteq \mcal{U} \cap \mcal{O}$ be a compact neighborhood of the origin with nonempty interior, and write $\mcal{N}_0 = \mcal{N}\setminus\left\{0 \right\}$.  Suppose $u(x)$ is a stationary control which is $C^1$ on $\mcal{N}_0$ and continuous on $\mcal{N}$.  Writing $F_u(x)=f\cp{x,u(x)}$, observe that we then similarly have 
    
    \vspace*{-0.5cm}\begin{align*}
    \restr{J_{F_u}}{x} &= \restr{\frac{\partial f}{\partial x}}{(x,u(x))}+\restr{\frac{\partial f}{\partial u}}{(x,u(x))}\restr{\frac{\partial u}{\partial x}}{x} = \restr{\frac{\partial f}{\partial u}}{(x,u(x))} \left ( \restr{\frac{\partial u}{\partial x}}{x} - Q\cp{x,u(x)} \right ).
    \end{align*}
    
    \vspace*{-0.25cm}Combining the above observations and setting $v=u(y)$, it follows that
    
    \vspace*{-0.5cm}\begin{align}
        \left (-\restr{J_{F[y,u(y)]}}{0}+\restr{J_{F_u}}{y}\right) &= \restr{\frac{\partial f}{\partial u}}{(y,u(y))}\left ( Q\cp{y,u(y)} - \restr{\frac{\partial w[y,u(y)]}{\partial x}}{0} \right )+\restr{\frac{\partial f}{\partial u}}{(y,u(y))} \left ( \restr{\frac{\partial u}{\partial x}}{y} - Q\cp{y,u(y)} \right ) \nonumber\\ 
        I+\restr{J_{F_u}}{y} &= \restr{\frac{\partial f}{\partial u}}{(y,u(y))}\left (\restr{\frac{\partial u}{\partial x}}{y} - \restr{\frac{\partial w[y,u(y)]}{\partial x}}{0}\right ). \label{auxeq}
    \end{align}
    
    \vspace*{-0.15cm}Therefore, if there exists a control $u(x)$ such that, for all $y \in \mcal{N}_0$,  
    
    \vspace*{-0.15cm}$$\restr{\frac{\partial f}{\partial u}}{(y,u(y))}\restr{\frac{\partial u}{\partial x}}{y} = \restr{\frac{\partial f}{\partial u}}{(y,u(y))}\restr{\frac{\partial w[y,u(y)]}{\partial x}}{0},$$
    
    \noindent it then follows that $\restr{J_{F_u}}{y}=\restr{J_{F[y,u(y)]}}{0} = -I$ for all $y \neq 0$ in $\mcal{N}_0$.  Moreover, observe that 
    
    \vspace*{-0.15cm}$$V(x)=\left \{ \begin{matrix} 
    \frac{1}{2}F_u(x)^*F_u(x) & x \neq 0 \\
    0 & x=0\end{matrix} \right.$$
    
    \noindent is then a Lyapunov function for $\dot{x}=F_u(x)$, as 
    
    \vspace*{-0.15cm}$$\langle \nabla V(x), F_u(x) \rangle = F_u(x)^*\restr{J_{F_u}}{x}F_u(x) = - \|F_u(x)\|^2.$$
    
    Finally, observe that for all $y \in \mcal{N}_0$, semitransversality yields a right inverse for $\restr{\frac{\partial f}{\partial u}}{(y,u(y))}$.  Thus, to yield a control $u(x)$ that satisfies
    
    $$\restr{\frac{\partial f}{\partial u}}{(y,u(y))}\restr{\frac{\partial u}{\partial x}}{y} = \restr{\frac{\partial f}{\partial u}}{(y,u(y))}\restr{\frac{\partial w[y,u(y)]}{\partial x}}{0}$$ 
    
    \noindent for all $y \in \mcal{N} \setminus \left \{ 0 \right \}$, we need to only solve
    
    \vspace*{-0.25cm}\begin{align*}
       \restr{\frac{\partial u}{\partial x}}{y} &= \left(\restr{\frac{\partial f}{\partial u}}{(y,u(y))}\right)^+\restr{\frac{\partial f}{\partial u}}{(y,u(y))}\restr{\frac{\partial w[y,u(y)]}{\partial x}}{0}.
    \end{align*}
\end{proof}

The lemma above establishes that the existence of a control satisfying some differentiability conditions and the equation \eqref{maineq} is sufficient to guarantee exponential stabilizability under the conditions of Theorem~\ref{ST,II}.  The next lemma will show that any such control satisfying \eqref{maineq} is a fixed point of a certain continuous self-mapping of a Banach space.

\begin{lemma}\label{supporting lemma 2}
    Suppose $f$ satisfies the conditions of Theorem~\ref{ST,II}. Let $\mcal{N} \subseteq \mcal{U}\cap \mcal{O}$ be a compact neighborhood of the origin with nonempty interior and let $\mcal{M}=\mcal{V} \cap \mcal{O}$.  Define the operator $T$ for $u \in C\cp{\mcal{N},\mcal{M}}$ by the action below, where the integral is understood component-wise:
    
    \vspace*{-0.25cm}$$[Tu](x) = \int_0^1 -\left(\restr{\frac{\partial f}{\partial u}}{\cp{sx,u(sx)}}\right)^+\left(I+\restr{\frac{\partial f}{\partial x}}{\cp{sx,u(sx)}} \right)x \es ds.$$
    
    \noindent Then, $T$ is a continuous self-mapping of the Banach space $C(\mcal{N},\mcal{M})$ which satisfies 
    $$T\cp{C(\mcal{N},\mcal{M})} \subseteq C^1\cp{\mcal{N}_0,\mcal{M} }\cap C\cp{\mcal{N},\mcal{M}}.$$  Moreover, $u(x) \in C^1\cp{\mcal{N}_0,\mcal{M} }\cap C\cp{\mcal{N},\mcal{M}}$ is stationary and satisfies \eqref{maineq} if and only if $u(x)$ is a fixed point of $T$.  
\end{lemma}
\begin{proof}
    Suppose $u(x) \in C^1\cp{\mcal{N}_0,\mcal{M} }\cap C\cp{\mcal{N},\mcal{M}}$.  We will first show that if $u(x)$ is stationary and satisfies \eqref{maineq}, $u(x)$ is a fixed point of $T$.  Using Serge Lang's vector-valued analogue of the mean value theorem \cite{lang}, observe that, for any $x\in \mcal{N}_0$ and any integer $k>0$, a solution $u(x)$ to \eqref{maineq} must satisfy the following, where $\gamma(s,k,x)=sx+\frac{1}{k}(1-s)x$ and the integral is understood component-wise:
    
    \vspace*{-0.5cm}\begin{align*}
        u(x)-u\cp{(1/k)x} &= \int_0^1 \restr{\frac{\partial u}{\partial x}}{\gamma(s,k,x)} ds \es \left(1-\frac{1}{k}\right) x \\
        &= \left(\frac{k-1}{k}\right)\int_0^1 \left(\restr{\frac{\partial f}{\partial u}}{\cp{\gamma(s,k,x),u(\gamma(s,k,x))}}\right)^+\restr{\frac{\partial f}{\partial u}}{\cp{\gamma(s,k,x),u(\gamma(s,k,x))}}\restr{\frac{\partial w[\gamma(s,k,x),u(\gamma(s,k,x))]}{\partial x}}{0} ds \es x \\
        &= -\left(\frac{k-1}{k}\right)\int_0^1 \left(\restr{\frac{\partial f}{\partial u}}{\cp{\gamma(s,k,x),u(\gamma(s,k,x))}}\right)^+\left(I+\restr{\frac{\partial f}{\partial x}}{\cp{\gamma(s,k,x),u(\gamma(s,k,x))}} \right) ds \es x.
    \end{align*}

    \noindent For notational convenience, let us set
    
    \vspace*{-0.15cm}$$U(x,y)=-\left(\restr{\frac{\partial f}{\partial u}}{\cp{x,u(x)}}\right)^+\left(I+\restr{\frac{\partial f}{\partial x}}{\cp{x,u(x)}} \right)y.$$
    
    The continuity of $U(x,y)$ on $\mcal{N}\times\mcal{N}$ is immediate from the continuity of $u$ and the continuous differentiability of $f$ on $\mcal{N}$.  So, for any fixed $s \in [0,1]$ and all $\epsilon>0$, there exits $\delta_\epsilon >0$ such that $\left\|U\cp{\gamma(s,k,x),x}-U(sx,x)\right\| < \epsilon$ whenever 
    
    \vspace*{-0.15cm}$$\|(\gamma(s,k,x),x)-(sx,x)\|=\|sx+(1/k)(1-s)x-sx\| = \left \| \cp{(1-s)/k}x \right \| < \delta_\epsilon.$$  
    
    \noindent Write $K_\mcal{N}=\sup{r > 0 :r\mbb{B}_{\R^n} \subseteq \mcal{N}}$ and $K_\mcal{M}=\sup{r > 0 :r\mbb{B}_{\R^m} \subseteq \mcal{M}}$.  Then, choosing $k \geq \delta_\epsilon/K_\mcal{N}$, it follows that $\| \frac{1-s}{k}x \| < \delta_\epsilon$ for any $x \in \mcal{N}$, and we have that $\left\|U\cp{\gamma(s,k,x),x}-U(sx,x)\right\| < \epsilon$ for $k>\delta_\epsilon/K_\mcal{N}$ and all $x \in \mcal{N}$.  Thus, $U(\gamma(s,k,x),x)$ converges to $U(sx,x)$ uniformly on $\mcal{N}$, and yields that any stationary solution $u(x) \in C^1\cp{\mcal{N}_0,\mcal{M} }\cap C\cp{\mcal{N},\mcal{M}}$ to \eqref{maineq} must satisfy
    
    \vspace*{-0.5cm}\begin{align*}
    u(x) &= \int_0^1 U(sx,x) ds + u(0)=[Tu](x).
    \end{align*}
    
    Since the continuity of $T$ is apparent, we will now show that $T$ is a self-mapping of the Banach space $C\cp{\mcal{N},\mcal{M}}$ satisfying $T\left ( C\cp{\mcal{N},\mcal{M}} \right ) \subseteq C^1\cp{\mcal{N}_0,\mcal{M} }\cap C\cp{\mcal{N},\mcal{M}}$. Without loss of generality, assume 
    
    \vspace*{-0.25cm}$$\txt{sup} \left \{ \left \|\left(\restr{\frac{\partial f}{\partial u}}{\cp{x,u}}\right)^+\left(I+\restr{\frac{\partial f}{\partial x}}{\cp{x,u}} \right)\right \| : (x,u) \in \mcal{O} \right \} \leq \frac{K_\mcal{M}}{K_\mcal{N}},$$ 
    
    \noindent as we could instead consider $cf(x,u)$ instead for some appropriate normalizing factor $c>0$ (or modify the norm in the above similarly) if this is not the case.  So, for $u\in C\cp{\mcal{N},\mcal{M}}$, it follows that
    
    \vspace*{-0.25cm}$$\sup{\|[Tu](x)\| : x \in \mcal{N}} \leq \sup{\frac{K_\mcal{M}}{K_\mcal{N}} \|x\| : x \in \mcal{N}} = K_\mcal{M}. $$ 
    
    Thus, we may conclude $T\cp{C(\mcal{N},\mcal{M})}\subseteq C(\mcal{N},\mcal{M})$. We show that $T\cp{C(\mcal{N},\mcal{M})} \subseteq C^1\cp{\mcal{N}_0,\mcal{M} }\cap C\cp{\mcal{N},\mcal{M}}$ by noting that for any $u \in C(\mcal{N},\mcal{M})$ and any $x\in \mcal{N}_0$, we have 
    
    \vspace*{-0.5cm}\begin{align*}
        \lim\limits_{\|h\|\rightarrow 0} \frac{\left \|[Tu](x+h)-[Tu](x)-U(x,h) \right \|}{\|h\|} &= \lim\limits_{\|h\|\rightarrow 0} \frac{\left \|\int_0^1 U\cp{s(x+h),x+h} ds -\int_0^1 U(sx,x) ds -U(x,h) \right \|}{\|h\|}  \\
        &=\lim\limits_{\|h\|\rightarrow 0}  \frac{\left \|\int_0^{x+h} \restr{\frac{\partial U}{\partial y}}{\cp{y,0}}dy - \int_0^x\restr{\frac{\partial U}{\partial y}}{\cp{y,0}} dy -U(x,h) \right \|}{\|h\|}  \\
        &=\lim\limits_{\|h\|\rightarrow 0}  \frac{\left \|\int_x^{x+h} \restr{\frac{\partial U}{\partial y}}{\cp{y,0}} dy -U(x,h) \right \|}{\|h\|}  \\
        &=\lim\limits_{\|h\|\rightarrow 0}  \frac{\left \| U(x,h)-U(x,h) \right \|}{\|h\|} = 0.
    \end{align*}
    
    Now, let us show that a fixed point of $T$ is a stationary solution to \eqref{maineq}.  To do so, first note that since $[Tu](0)=0$, any fixed point of $T$ is necessarily stationary.  Then, let us note that the manipulation above showing differentiability yields the following for any $y \in \mcal{N}_0$: 
    $$\restr{\frac{\partial [Tu]}{\partial x}}{y}=-\left(\restr{\frac{\partial f}{\partial u}}{\cp{y,u(y)}}\right)^+\left(I+\restr{\frac{\partial f}{\partial x}}{\cp{y,u(y)}} \right).$$  
    Now, note that \eqref{auxeq} gives that 
    \begin{align*}
        I+\restr{J_{F_u}}{y} &= \restr{\frac{\partial f}{\partial u}}{(y,u(y))}\left (\restr{\frac{\partial u}{\partial x}}{y} - \restr{\frac{\partial w[y,u(y)]}{\partial x}}{0}\right ) \\
       I+\restr{\frac{\partial f}{\partial x}}{(y,u(y))}+\restr{\frac{\partial f}{\partial u}}{(y,u(y))}\restr{\frac{\partial u}{\partial x}}{y} &= \restr{\frac{\partial f}{\partial u}}{(y,u(y))}\restr{\frac{\partial u}{\partial x}}{y} - \restr{\frac{\partial f}{\partial u}}{(y,u(y))}\restr{\frac{\partial w[y,u(y)]}{\partial x}}{0} \\
        I+\restr{\frac{\partial f}{\partial x}}{(y,u(y))}&= - \restr{\frac{\partial f}{\partial u}}{(y,u(y))}\restr{\frac{\partial w[y,u(y)]}{\partial x}}{0} \\
    \end{align*}
    So, if $u$ is a fixed point of $T$, then $u$ must lie in $C^1\cp{\mcal{N}_0,\mcal{M} }\cap C\cp{\mcal{N},\mcal{M}}$ and correspondingly satisfy the following for all $y\in \mcal{N}_0$:
    \begin{align*}
        \restr{\frac{\partial u}{\partial x}}{y} &= \restr{\frac{\partial [Tu]}{\partial x}}{y} \\
        &= -\left(\restr{\frac{\partial f}{\partial u}}{\cp{y,u(y)}}\right)^+\left(I+\restr{\frac{\partial f}{\partial x}}{\cp{y,u(y)}} \right)\\
        &=\left(\restr{\frac{\partial f}{\partial u}}{(y,u(y))}\right)^+\restr{\frac{\partial f}{\partial u}}{(y,u(y))}\restr{\frac{\partial w[y,u(y)]}{\partial x}}{0}.
    \end{align*}
\end{proof}

\noindent The lemma above establishes that solutions to \eqref{maineq} which satisfy the conditions of Lemma~\ref{supporting lemma 1} are fixed points of a continuous self-mapping of a Banach space and vice versa.  Using these results together, we can finally prove Theorem~\ref{ST,II}
\begin{proof}
    By Lemma~\ref{supporting lemma 2}, a fixed point of $T$ yields a solution $u \in C^1\cp{\mcal{N}_0,\mcal{M} }\cap C\cp{\mcal{N},\mcal{M}}$ to \eqref{maineq}.  By Lemma~\ref{supporting lemma 1}, a solution $u \in C^1\cp{\mcal{N}_0,\mcal{M} }\cap C\cp{\mcal{N},\mcal{M}}$ to \eqref{maineq} gives our desired result.  So, we need only show that $T$ has a fixed point.  
    
    \bs
    
    To that end, denote by $H^{k}(\Omega)$ the usual Sobolev space of functions $q \in L^2(\Omega)$ with derivatives $D^\alpha q \in L^2(\Omega)$ for each multi-index $\alpha$ with $|\alpha|\leq k$. Write $X=C^1\cp{\mcal{N}_0,\mcal{M} }\cap C\cp{\mcal{N},\mcal{M}}$ and for $r > 0$, set 
    
    \vspace*{-0.25cm}$$\mcal{B}_X(r) = \left \{ u \in X : u_i \in H^{1}(\mcal{N}), \max_{i\in[1,n]}\|u_i\|_{H^{1}(\mcal{N})}\leq r \right \}.$$
    
    To determine the closure of $\mcal{B}_X(r)$ in $C\cp{\mcal{N},\mcal{M}}$, note that the closure of $\prod_{i=1}^n r\mbb{B}_{H^1(\mcal{N})}$ in $\prod_{i=1}^n r\mbb{B}_{H^0(\mcal{N})}$ is simply 
    
    $$cl\left ( \prod_{i=1}^n r\mbb{B}_{H^1(\mcal{N})} \right ) = \prod_{i=1}^n r\mbb{B}_{H^0(\mcal{N})},$$ 
    
    \noindent where the product is understood as an $n$-fold Cartesian product and the closure is understood in the resulting product topology.  By the Rellich-Kondrachov theorem, this closure is compact, so observe that $X \cap \prod_{i=1}^n r\mbb{B}_{H^1(\mcal{N})}$ is dense in $\prod_{i=1}^n r\mbb{B}_{H^1(\mcal{N})}$, yielding  
    
    \vspace*{-0.25cm}$$cl\cp{\mcal{B}_X(r)} = \prod_{i=1}^n r\mbb{B}_{H^0(\mcal{N})} \cap C\cp{\mcal{N},\mcal{M}}.$$
    
    \noindent Hence, $cl\cp{\mcal{B}_X(r)}$ is compact in $C\cp{\mcal{N},\mcal{M}}$.  Further, the inclusion $T\cp{C(\mcal{N},\mcal{M})} \subseteq X$ yields that $T\cp{cl\cp{\mcal{B}_X(r)}} \subseteq \mcal{B}_X(r)$, and it is routine to verify that $s(u(x)+(1-s)v(x) \in cl\cp{\mcal{B}_X(r)}$ for all $s \in [0,1]$ and all $u,v\in cl\cp{\mcal{B}_X(r)}$.  Therefore, $cl\cp{\mcal{B}_X(r)}$  is convex and compact, and $T$ is a continuous mapping which sends the nonempty, compact, convex subset of the Banach space $C\cp{\mcal{N},\mcal{M}}$ into itself.  By the Tychonoff fixed point theorem, the result then follows.
\end{proof}

Thus we reveal that \textit{linear openness} and \textit{semitransversality} in a \textit{punctured neighborhood} of equilibrium point is sufficient to account for whatever deficiencies the linearization may present.  The obtained result tells us that, by strengthening the openness property in Brockett's necessary condition (Theorem \ref{brocketts}) to the one with a {\em linear rate}, we achieve the local exponential stabilizability of smooth control systems by means of continuous stationary feedback laws {\em with only one additional} (minimally restrictive) {\em assumption}, in contrast to the previous result of \cite{GJKM} in Theorem~\ref{GJKM 3.2}.  


\section{Concluding Remarks}\label{sec5}


This paper reveals, by using the machinery of {\em variational analysis}, that a \textit{linear openness} extension of Brockett's well-known \textit{openness necessary} condition is \textit{sufficient} for a system to be locally exponentially stabilizable under the action of a \textit{composition operator} induced by a local {\em diffeomorphism}. This largely extends recent results of \cite{GJKM} obtained under restrictive spectral assumptions. Further, it is shown in this paper that the \textit{local exponential stabilizability} of the system under the action of this composition operator can be related to both \textit{local exponential stabilizability} and \textit{local asymptotic stabilizability} of the \textit{original system} via the novel \textit{semitransversality condition}, and that this condition proves not only sufficient \textit{point-wise} at an equilibrium, but \textit{locally} near an equilibrium.\vspace*{0.05in}

We see now some directions for further developments of the obtained results:\vspace*{0.05in}
\begin{enumerate}[(i)]
\item Extending the established characterizations of local asymptotic and exponential stabilizability of linearly open ODE systems to the case of {\em nonsmooth} vector fields. As mentioned in Section~\ref{sec2}, variational analysis achieves complete characterizations of linear openness for general nonsmooth mappings, and the challenge is to implement them for the study of feedback stabilizability.\vspace*{0.05in}
\item Theorem~\ref{QTN,I} reveals that the existence of a feedback law locally exponentially stabilizing a certain modified class of systems is equivalent to the existence of a stabilizing composition operator for \eqref{sys}.  While this may be theoretically of interest, it also suggests possible concrete, practical applications to Lyapunov theoretic approaches to stabilizability.  Namely, in the context of Theorem~\ref{QTN,I}, suppose $v \in \R^k$ for any $k \geq 0$. Write $w=\ba{c}{x \\ v}$, and suppose $g(x,v)$ is any system satisfying $covg(0,0)>0$. Then, setting $F(x,w) = f\cp{x,u^*(x)}-cg(w)$ for some continuously differentiable feedback law $u^*(x)$ and some constant $c>0$, the locally exponential stability of $\dot{x}=F(x,w(x))$ then yields the existence of a Lyapunov function $V_F$ satisfying $\left \langle \nabla V_F(x), F\cp{x,w(x)} \right \rangle  < 0 $ for all $x$ in a neighborhood of the origin.  But this then produces, by the construction of $F$, that $\left \langle \nabla V_F(x), f\cp{x,u(x)} \right \rangle < \left \langle \nabla V_F(x), cg\cp{w(x)} \right \rangle$.  So it then follows that if, for whatever reason, it can be shown that $\dot{x}=cg\cp{w(x)}$ is locally asymptotically stable and a Lyapunov function $V_g$ for $cg\cp{w(x)}$ can be constructed to satisfy $\left \langle \nabla V_F(x), cg\cp{w(x)} \right \rangle < \left \langle \nabla V_F(x), cg\cp{w(x)} \right \rangle$ for all $x \neq 0$ in some neighborhood of the origin, it then follows that $u(x)$ is a stabilizing control for $\dot{x}=f(x,u)$.  While this remains an observation at present, the fact that $\dot{x}=g(x,v)$ need only be stabilized \textit{by a composition operator} $w(x)$ instead of the usual, more restrictive feedback law formulation suggests such an approach may be worthy of consideration.  Particularly, since exponential stabilizability by $C^1$ feedback laws is already well understood (see \cite{zab}), such an approach may provide a route towards characterizing the possible Lyapunov functions for a given system.  However, further significant progress in this direction requires an improvement of the row-rank condition  on $B_f$ imposed in Lemma~\ref{QT,LII} to characterize semitranversality. \vspace*{0.05in}
\end{enumerate}


\bbs

\begin{addmargin}[1em]{1em}
\small{
\noindent \textothersc{BC}: \hspace*{0.2cm} \textothersc{Department of Mathematics and Statistics, \newline \hspace*{0.85cm} University of Wyoming, Laramie, WY 82071}

\hspace*{0.22cm} \textit{E-mail Address}: \textbf{bchris19@uwyo.edu}

\bs

\noindent \textothersc{FJ}: \hspace*{0.3cm} \textothersc{Department of Mathematics and Statistics, \newline \hspace*{0.8cm} University of Wyoming, Laramie, WY 82071}

\hspace*{0.22cm} \textit{E-mail Address}: \textbf{fjafari@uwyo.edu}

\bs

\noindent \textothersc{BM}: \es \textothersc{Department of Mathematics, \newline \hspace*{0.8cm} Wayne State University, Detroit, MI 48202}

\hspace*{0.22cm} \textit{E-mail Address}: \textbf{boris@math.wayne.edu}
}
\end{addmargin}


\begin{thebibliography}{99}
\bibitem{andreini} A. Andreini, A. Bacciotti and  G. Stefani, \textit{Global Stabilizability of Homogeneous Vector Fields of Odd Degree}, Systems \& Control Letters \textbf{10} (1988), 251-256.
\bibitem{artstein} Z. Artstein, \textit{Stabilization with Relaxed Controls}, Nonlinear Anal. \textbf{7} (1983), 1163-1173.
\bibitem{brockett} R. W. Brockett, \textit{Asymptotic Stability and Feedback Stabilization}, In \textit{Differential Geometric Control Theory} (R. W. Brockett et al., eds.), pp. 181?191, Birkh\"auser, Boston, MA, 1983.
\bibitem{byrnes} C. I. Byrnes, \textit{On Brockett?s Necessary Condition for Stabilizability and the Topology of Lyapunov Functions on $\R^N$}, Commun. Inf. Syst. \textbf{8} (2008), 333-352.
\bibitem{coron1} J. M. Coron, \textit{Control and Nonlinearity}, American Mathematical Society, Providence, RI, 2007.
\bibitem{coron2} J. M. Coron, \textit{A Necessary Condition for Feedback Stabilization}, Syst. Control Lett. \textbf{14} (1990), 227-232.
\bibitem{milyutin} A. V. Dmitruk, A. A. Milyutin and N. P. Osmolovskii, \textit{Lyusternik?s Theorem and the Theory of Extrema}, Russian Math. Surveys \textbf{35} (1980), 11-51.
\bibitem{GJKM} R. Gupta, F. Jafari, R. J. Kipka and B. S. Mordukhovich.  \textit{Linear Openness and Feedback Stabilization of Nonlinear Control Systems}, Discrete \& Continuous Dynamical Systems \textbf{11} (2018), 1103-1119.
\bibitem{hermes}  H. Hermes, \textit{Asymptotically Stabilizing Feedback Controls and the Nonlinear Regulator Problem}, SIAM J. Control Optim. \textbf{29} (1991), 185-196.
\bibitem{ioffe} A. D. Ioffe, \textit{Metric Regularity -- A Survey, Part 1: Theory and Applications}, Journal of the Australian Mathematical Society \textbf{101} (2016), 188-243.
\bibitem{ishikawa} M. Ishikawa, M. Sampei, \textit{On Equilibria Set and Feedback Stabilizability of Nonlinear Control Systems}, IFAC Proceedings \textbf{31} 1998, 609-614
\bibitem{lang} S. Lang, \textit{Analysis I}, Addison-Wesley Publishing Company (1986)
\bibitem{mordukhovich1}  B. S. Mordukhovich, \textit{Complete Characterization of Openness, Metric Regularity, and Lipschitzian Properties of Multifunctions}, Trans. Amer. Math. Soc. \textbf{340} (1993), 1-35.
\bibitem{mordukhovich2} B. S. Mordukhovich, \textit{Variational Analysis and Generalized Differentiation, I: Basic Theory}, Springer, Berlin, 2006.
\bibitem{mordukhovich3} B. S. Mordukhovich, \textit{Variational Analysis and Applications}, Springer, Cham, Switzerland, 2018.
\bibitem{rockafellar} R. T. Rockafellar and R. J-B. Wets, \textit{Variational Analysis}, Springer, Berlin, 1998.
\bibitem{sontag2} E. D. Sontag, \textit{Mathematical Control Theory: Deterministic Finite Dimensional Systems}, Springer, New York, 1998.
\bibitem{sontag3} E. D. Sontag, \textit{Stability and Stabilization: Discontinuities and the Effect of Disturbances}. In \textit{Nonlinear Analysis, Differential Equations and Control} (F. H. Clarke and R. J. Stern, eds.), pp. 551-598, Kluwer, Dordrecht, The Netherlands, 1999.
\bibitem{sontag1} H. J. Sussmann, E. D. Sontag and D. Y. Yang, \textit{A General Result on the Stabilization of Linear Systems using Bounded Controls}, IEEE Trans. Automat. Control \textbf{39} (1994), 2411-2425.
\bibitem{zab} J. Zabczyk, \textit{Some Comments on Stabilizability},  J. Applied Math and Optimization \textbf{19} (1989), 1-9
\end{thebibliography}
\end{document}